\documentclass[final]{siamltex}%
\usepackage{amsmath,amssymb}
\usepackage{amsfonts}
\usepackage{exscale}
\usepackage{graphicx}
\usepackage{multirow}
\usepackage{url}
\usepackage[pdftex,
pdfstartview=FitH,bookmarks=true,bookmarksnumbered=true,bookmarksopen=true,hypertexnames=false,breaklinks=true,
colorlinks=true,  linkcolor=blue,anchorcolor=blue,
citecolor=blue,filecolor=blue,  menucolor=blue]{hyperref}
\usepackage{amsmath}
\usepackage{amssymb}
\usepackage{version}%
\providecommand{\U}[1]{\protect\rule{.1in}{.1in}}

\newtheorem{algorithm}[theorem]{Algorithm}
\newtheorem{assumption}[theorem]{Assumption}

\newtheorem{remark}[theorem]{Remark}
\makeatletter
\def\@cite#1#2{{\rm [}{{\rm#1}\if@tempswa , #2\fi}{\rm ]}}
\makeatother
\begin{document}
\title{On adaptive {BDDC} for the flow \\
in heterogeneous porous media \thanks
{Supported by the U.S. National Science Foundation under grant DMS1521563.}}
\author{Bed\v{r}ich Soused\'{\i}k \thanks
{Department of Mathematics and Statistics,
University of Maryland, Baltimore County, 1000 Hilltop Circle, Baltimore, MD~21250, USA,
(\texttt{sousedik@umbc.edu}) }}
\maketitle\begin{abstract}
We study a method based on Balancing Domain Decomposition by Constraints
({BDDC)}
for a numerical solution of a single-phase flow in heterogenous porous media. The
method solves for both flux and pressure variables.
The fluxes are resolved in three steps: the coarse solve is followed by subdomain solves
and last we look for a divergence-free flux correction and pressures
using conjugate gradients with the BDDC preconditioner.
Our main contribution is an application of the adaptive algorithm for selection of flux constraints.
Performance of the method is
illustrated on the benchmark problem from the 10th {SPE} {C}omparative
{S}olution {P}roject\ ({SPE~10}).
Numerical experiments in both 2D and 3D demonstrate that the first two steps of the method exhibit some numerical upscaling properties, and the adaptive preconditioner in the last step
allows a significant decrease in number of iterations of conjugate gradients at a small additional cost.
\end{abstract}
\begin{AMS}
{65F08, 65F10, 65M55, 65N55}
\end{AMS}
\pagestyle{myheadings}
\thispagestyle{plain}
\markboth{BED\v{R}ICH SOUSED\'{I}K}{ON ADAPTIVE BDDC FOR FLOW IN POROUS MEDIA}%

\section{Introduction}

\label{sec:introduction}The {B}alancing {D}omain {D}ecomposition by
{C}onstraints (BDDC), proposed independently by Cros~\cite{Cros-2003-PSC},
Dohrmann~\cite{Dohrmann-2003-PSC}, and Fragakis and
Papadrakakis~\cite{Fragakis-2003-MHP}, is one of the most popular methods of
iterative substructuring. The method was developed as a preconditioner for the
solution of systems of linear equation obtained by finite element
discretizations of elliptic problems, and it has been originally derived as a
primal counterpart of the Finite Element Tearing and Interconnecting - Dual,
Primal (FETI-DP) method by Farhat et
al.~\cite{Farhat-2001-FDP,Farhat-2000-SDP}. Over the years the BDDC has been
extended to other types of problems, for example to the nearly incompressible
elasticity by Dohrmann~\cite{Dohrmann-2004-SSP}, the Stokes problem by Li and
Widlund~\cite{Li-2006-BAI}, or advection-diffusion problems by Tu and
Li~\cite{Li-2009-CAB,Tu-2008-BDD}. It is also relatively straightforward to
extend the BDDC\ into multiple levels, as noted by
Dohrmann~\cite{Dohrmann-2003-PSC}. The three-level methods were developed in
two and three dimensions by Tu~\cite{Tu-2007-TBT3D,Tu-2007-TBT}, and Mandel et
al.~\cite{Mandel-2008-MMB} extended the algorithm into a multilevel method
within a more general multispace BDDC\ setting. Another class of problems,
important in the context of this paper, is the flow in porous media based on
\emph{mixed} and \emph{mixed-hybrid} finite element discretizations. Probably
the first domain decomposition methods of this class were proposed by
Glowinski and Wheeler~\cite{Glowinski-1988-DDM}. Their Method~II was
preconditioned using BDD by Cowsar et al.~\cite{Cowsar-1995-BDD}, using BDDC
by Tu~\cite{Tu-2007-BAF}, and \v{S}\'{\i}stek et al.~\cite{Sistek-2015-BMF}
extended this methodology to flow in porous media with combined mesh
dimensions. This approach is regarded as \emph{hybrid} because the method
iterates on a system of \emph{dual} variables (as Lagrange multipliers)
enforcing the continuity of flux variables across the substructure interfaces.
An alternative strategy, retaining the original \emph{primal} variables was
proposed by Tu~\cite{Tu-2005-BAM,Tu-2011-TBA}, who combined the BDDC
preconditioner with an earlier algorithmic framework developed by Ewing and
Wang~\cite{Ewing-1992-ASA}, cf. also Mathew~\cite{Mathew-1993-SAIa}, which
allows to solve the saddle-point problem obtained from \emph{mixed} finite
element discretization by conjugate gradients. The Nested BDDC by
Soused\'{\i}k~\cite{Sousedik-2013-NBS} provided a multilevel extension by
applying the framework from~\cite{Tu-2005-BAM} recursively. Most recently,
Zampini and Tu~\cite{Zampini-2017-MBD} presented another approach to
multilevel BDDC\ including adaptive coarse space construction, which relies on
a special, so-called, deluxe scaling.

There are two main ingredients of a BDDC\ method: a coarse space, which is
defined by \emph{constraints} on the values of degrees of freedom, and a
scaling (averaging) operator, which provides a mapping between the solution
space and the space in which the solves in the preconditioner are performed.
The algorithm for adaptive selection of constraints for both methods the BDDC
and FETI-DP\ was originally proposed by Mandel and
Soused\'{\i}k~\cite{Mandel-2007-ASF}.\ The algorithm was later generalized in
a joint work with \v{S}\'{\i}stek~\cite{Mandel-2012-ABT} into three spatial
dimensions and implemented for the BDDC using an approach inspired by a
partial subassembly and a change of variables by Li and
Widlund~\cite{Li-2006-FBB}. Finally, we also reformulated the algorithm to
treat the coarse space explicitly~\cite{Sousedik-2013-AMB}. We note that there
are many other approaches to the adaptive construction of the coarse spaces in
BDDC, see~\cite{Pechstein-2017-UFA} and the references therein, as well as for
BDD see, e.g.,~\cite{Spillane-2013-ASC}. There have been several scalings
studied in the literature. In the \emph{multiplicity} scaling, the weights are
chosen proportionally to the number of subdomains sharing a given degree of
freedom, and it is regarded as not robust for coefficient jumps. The $\rho
$-\emph{scaling }leads to robustness, but it relies on knowledge of the
problem coefficients~\cite{Klawonn-2008-AFA}. The \emph{stiffness} scaling is based on the diagonal of
the stiffness matrix, but in some cases with irregular meshes it may lead to
high condition numbers~\cite{Pechstein-2013-FBE,Pechstein-2011-AFM}. All these
scalings involve diagonal matrices. Finally, the \emph{deluxe} scaling
introduced in~\cite{Dohrmann-2013-SRT} uses dense matrices, which are computed
from inverses of localized Schur complements. It has been observed to be quite
robust~\cite{Oh-2018-BAD,Zampini-2017-MBD} but also computationally intensive.

In this paper, we build on the primal strategy. The starting point is the
two-level algorithm from~\cite{Sousedik-2013-NBS}, which we combine with
adaptive selection of constraints following~\cite{Sousedik-2013-AMB} and apply
it to flow in heterogenous porous media. To this end, we use a reservoir from
the 10th {SPE} {C}omparative {S}olution {P}roject\ ({SPE~10}) cf.,
e.g.,~\cite{Aarnes-2007-INF,Christie-2001-SPE10} as the benchmark problem. The
BDDC method from~\cite{Sousedik-2013-NBS} solves for both flux and pressure
variables. The fluxes are resolved in three steps: the coarse solve is
followed by mutually independent subdomain solves, and last we look for a
divergence-free flux correction and pressure using conjugate gradients (CG)
with the BDDC preconditioner. The coarse solve in the first step is exactly
the same as the coarse solve used in the BDDC preconditioner in the step
three. It is assumed that the initial constraints preserve the iterates in a
\emph{balanced} subspace, in which the preconditioned operator is positive
definite. Our goal here is to adapt the method to flow in realistic
reservoirs, characterized by highly heterogeneous permeability coefficients in
as simple way as possible. In particular, we translate the ideas used for
elliptic problems {in}~\cite{Sousedik-2013-AMB} to mixed formulations of flow
in porous media discretized by the lowest-order Raviart-Thomas finite elements
(RT0). The main component of the extension\ is the use of additional
adaptive\ flux coarse basis functions. The starting point is the condition
number bound formulated as a generalized eigenvalue problem, which is replaced
by a number of local eigenvalue problems formulated for pairs of adjacent
subdomains, and the eigenvectors, corresponding to the eigenvalues larger than
a target condition number are used to construct the additional flux coarse
basis functions. We note that from this perspective our method can be viewed
as a way of numerical upscaling via the coarse basis functions known from the
BDDC. Unlike~\cite{Zampini-2017-MBD} we do not use a change of basis and
partial assembly of operators, and we also illustrate that for this problem
the multiplicity scaling in combination with the adaptive algorithm and a
simple diagonal rescaling of the pressure block in the setup of the problem is
sufficient to construct a robust algorithm. Numerical experiments in both 2D
and 3D demonstrate that the first two steps of the method exhibit some
numerical upscaling properties, and the convergence rate of conjugate
gradients in the last step can be estimated a priori in the setup of the
adaptive algorithm.

The paper is organized as follows. In Section~\ref{sec:model} we introduce the
model problem, in Section~\ref{sec:two-level}\ we recall the BDDC\ method and
the preconditioner, in Section~\ref{sec:adaptive}\ we formulate the algorithm
for adaptive selection of the flux constraints, in
Section~\ref{sec:implementation}\ we discuss some details of implementation,
in Section~\ref{sec:numerical}\ we present results of numerical experiments,
and finallly in Section~\ref{sec:conclusion} we summarize and conclude our work.

For convenience, we identify finite element functions with the vectors of
their coefficients in the corresponding finite element basis. These
coefficients are also called \emph{variables} or \emph{degrees of freedom}. At
a few places we will also identify linear operators with their matrices, in
bases that will be clear from the context. For a symmetric positive definite
bilinear form $a$, we will denote the energy norm by $\left\Vert u\right\Vert
_{a}=\sqrt{a\left(  u,u\right)  }$.

\section{Model problem}

\label{sec:model}Let $\Omega$ be a bounded domain in $%
\mathbb{R}
^{n}$, where $n=2$ or $3$. We~would like to find the solution of the following
mixed problem, which combines the Darcy's law
relating flux~$\mathbf{u}$ and pressure~$p$, and the
equation of continuity,
\begin{align}
k^{-1}\mathbf{u}+\nabla p  &  =0\quad\text{in }\Omega,\label{eq:problem-1}\\
\nabla\cdot\mathbf{u}  &  =f_{\Omega}\quad\text{in }\Omega,
\label{eq:problem-2}\\
p  &  =p_{N},\quad\text{on }\partial\Gamma_{N},\label{eq:problem-3}\\
\mathbf{u}\cdot\mathbf{n}  &  =g_{E}\quad\text{on }\partial\Gamma_{E},
\label{eq:problem-4}%
\end{align}
where $\partial\Omega=\overline{\Gamma}_{E}\cup\overline{\Gamma}_{N}$, and
$\mathbf{n}$\ denotes the unit outward normal of $\Omega$. The coefficient
$k=k_{p}/\mu$, where $k_{p}$ is the permeability of the porous medium
and\ $\mu$\ is the viscosity of the fluid. For simplicity, we will set $\mu=1$
and so$~k=k_{p}$. Without loss of generality we will also assume that
$\Gamma_{N}=\emptyset$, which requires a compatibility condition%
\begin{equation}
- \int_{\Omega}f_{\Omega}\,dx+\int_{\partial\Omega}g_{E}\,ds=0,
\label{eq:compatibility}%
\end{equation}
and the pressure $p$ will be uniquely determined up to an additive constant.
We will further assume that $g_{E}=0$. These assumptions motivate the
definition of a space%
\[
\mathbf{H}_{0}(\Omega;\operatorname{div})=\left\{  \mathbf{v:v}\in
L^{2}(\Omega);\nabla\cdot\mathbf{v}\in L^{2}(\Omega)\quad\text{and}%
\quad\mathbf{v}\cdot\mathbf{n}=0\;\text{on }\partial\Omega\right\}  ,
\]
equipped with the norm
\[
\left\Vert \mathbf{v}\right\Vert _{\mathbf{H}_{0}(\Omega;\operatorname{div}%
)}^{2}=\left\Vert \mathbf{v}\right\Vert _{L^{2}(\Omega)}^{2}+H_{\Omega}%
^{2}\left\Vert \nabla\cdot\mathbf{v}\right\Vert _{L^{2}(\Omega)}^{2},
\]
where $H_{\Omega}$ is the characteristic size of $\Omega$, and the definition
of a space%
\[
L_{0}^{2}\left(  \Omega\right)  =\left\{  q:q\in L^{2}\left(  \Omega\right)
\quad\text{and}\quad\int_{\Omega}qdx=0\right\}  .
\]
The weak form of the problem we wish to solve, is
\begin{align}
\int_{\Omega}k^{-1}\mathbf{u}\cdot\mathbf{v}\,dx-\int_{\Omega}p\left(
\nabla\cdot\mathbf{v}\right)  \,dx  &  =0,\quad\forall\mathbf{v}\in
\mathbf{H}_{0}(\Omega;\operatorname{div}),\label{eq:problem-weak-1}\\
-\int_{\Omega}\left(  \nabla\cdot\mathbf{u}\right)  q\,dx  &  =-\int_{\Omega
}f_{\Omega}q\,dx,\quad\forall q\in L_{0}^{2}\left(  \Omega\right)  .
\label{eq:problem-weak-2}%
\end{align}
We refer, e.g., to the monographs~\cite{Brezzi-1991-MHF,Toselli-2005-DDM} for
additional details and discussion.

Next, let $U$ be the lowest-order Raviart-Thomas (RT0)
finite element space with a zero normal component on $\partial\Omega$ and let
$Q$ be a space of piecewise constant finite element basis functions with a
zero mean on~$\Omega$.
These two spaces, defined on the triangulation$~\mathcal{T}_{h}$ of$~\Omega$,
where $h$\ denotes the mesh size, are finite dimensional subspaces of
$\mathbf{H}_{0}(\Omega;\operatorname{div})$ and $L_{0}^{2}(\Omega)$,
respectively, and they satisfy a uniform inf-sup condition,
see~\cite{Brezzi-1991-MHF}. Let us define the bilinear forms and the
right-hand side by%
\begin{align}
a\left(  u,v\right)   &  =\int_{\Omega}k^{-1}\mathbf{u}\cdot\mathbf{v}%
\,dx,\label{eq:a}\\
b\left(  u,q\right)   &  =-\int_{\Omega}\left(  \nabla\cdot\mathbf{u}\right)
q\,dx,\label{eq:b}\\
\left\langle f,q\right\rangle  &  =-\int_{\Omega}f_{\Omega}q\,dx.
\label{eq:rhs}%
\end{align}

In the mixed finite element approximation of problem (\ref{eq:problem-weak-1}%
)--(\ref{eq:problem-weak-2}), we would like to find a pair of fluxes and
pressures $\left(  u,p\right)  \in\left(  U,Q\right)  $ such that
\begin{align}
a\left(  u,v\right)  +b\left(  v,p\right)   &  =0,\qquad\forall v\in
U,\label{eq:variational-1}\\
b\left(  u,q\right)   &  =\left\langle f,q\right\rangle ,\qquad\forall q\in Q.
\label{eq:variational-2}%
\end{align}
We note that $Q$ is a finite-dimensional subspace of $L_{0}^{2}\left(
\Omega\right)  $ and therefore the unique solvability of the mixed problem
(\ref{eq:variational-1})--(\ref{eq:variational-2}) is guaranteed.

In the next section, we will describe the components of the two-level Nested
{BDDC}, which allows an efficient iterative solution of
problem~(\ref{eq:variational-1})--(\ref{eq:variational-2}).

\section{The BDDC method}

\label{sec:two-level}Let us consider a decomposition of$~\Omega$ into a set of
nonoverlapping subdomains $\Omega^{i}$, $i=1,\ldots,N,$ also called
substructures, forming a quasi-uniform triangulation of $\Omega$ and denote
the characteristic subdomain size by~$H$.
Each substructure is a union of finite elements with a matching discretization across the substructure
interfaces. Let $\Gamma^{i}=\partial\Omega^{i}\backslash\partial\Omega$\ be
the set of boundary degrees of freedom of a~substructure$~\Omega^{i}$\ shared
with another substructure$~\Omega^{j}$, $j\neq i$, and define the interface by
$\Gamma=\cup_{i=1}^{N}\Gamma^{i}$. Let us define a \emph{face} as an
intersection $\Gamma^{ij}=\Gamma^{i}\cap\Gamma^{j}$, $i\neq j$ and let us
denote by $\mathcal{F}$ the set of all faces between substructures. Note that
with respect to the~RT0 discretization we define only
\emph{faces}, but no \emph{corners} (nor \emph{edges} in 3D) known from other
types of substructuring.

We will solve problems similar to (\ref{eq:variational-1}%
)--(\ref{eq:variational-2}) on each substructure. As we have noted, such
problems determine the pressure uniquely up to a constant,
so we consider the decomposition of the pressure space
\begin{equation}
Q=Q_{0}\oplus Q_{I},\qquad Q_{I}=Q^{1}\times\cdots Q^{N},
\label{eq:Q_decomposition}%
\end{equation}
where $Q_{0}$ consists of functions that are constant in each subdomain and
have a zero average over the whole domain$~\Omega$, and the product
space$\ Q_{I}$ consists of functions that have zero weighted average over one
subdomain at a time. That is,
\begin{equation}
\int_{\Omega}q_{0}\,dx=0,\quad\forall q_{0}\in Q_{0}\text{\qquad and\qquad
}\int_{\Omega^{i}}q^{i}\,dx=0,\quad\forall q^{i}\in Q^{i},\;i=1,\dots,N.
\label{eq:Q_decomposition-int}%
\end{equation}

Next, let $W^{i}$ be the space of flux finite element functions on a
substructure $\Omega^{i}$ such that all of their degrees of freedom on
$\partial\Omega^{i}\cap\partial\Omega$ are zero, and let%
\[
W=W^{1}\times\dots\times W^{N}.
\]
Hence $U\subset W$ can be viewed as the subspace of
flux functions from~$W$ such that $u \cdot \mathbf{n}$ is continuous across substructure interfaces.
Define $U_{I}\subset U$ as the subspace of
flux functions such that $u \cdot \mathbf{n}$ is zero on the interface~$\Gamma$,
i.e., the space of \textquotedblleft interior\textquotedblright\ flux
functions, and let us also define a mapping
$P:w\in W\longmapsto u_{I}\in U_{I}$ such that%
\begin{align*}
a\left(  u_{I},v_{I}\right)  +b\left(  v_{I},p_{I}\right)   &  =a\left(
w,v_{I}\right)  ,\quad\forall v_{I}\in U_{I},\\
b\left(  u_{I},q_{I}\right)   &  =b\left(  w,q_{I}\right)  ,\quad\forall
q_{I}\in Q_{I}.
\end{align*}
Functions from $\left(  I-P\right)  W$
will be called Stokes harmonic, cf.~\cite[Section~9.4.2]{Toselli-2005-DDM}.

Let $\widehat{W}$ be the space of Stokes harmonic functions that are
continuous across substructure interfaces, and such that
\begin{equation}
U=\widehat{W}\oplus U_{I},\qquad\widehat{W}\perp_{a}U_{I}.
\label{eq:discrete-harm}%
\end{equation}
We note that from the divergence theorem, for all $u_{I}\in U_{I}$ and
$q_{0}\in Q_{0}$, we obtain%
\[
b\left(  u_{I},q_{0}\right)  =-\int_{\Omega}\left(  \nabla\cdot u_{I}\right)
q_{0}dx=0.
\]

The BDDC is a two-level method characterized by a selection of certain
\emph{coarse degrees of freedom}. In the present setting these will be flux
averages over faces shared by a pair of substructures at a time and pressure
averages over each substructure. Let us denote by $\widetilde{W}\subset\left(
I-P\right)  W$ the subspace of Stokes harmonic functions such that their flux
coarse degrees of freedom on adjacent substructures coincide; for this reason
we will use the terms coarse degrees freedom and \emph{constraints}
interchangeably. Specifically, we define a zero-net flux constraint for a
face$~\Gamma^{ij}$ as%
\begin{equation}
\int_{\Gamma^{ij}}\left(  w^{i}-w^{j}\right)  \cdot\mathbf{n}^{i}%
\mathbf{\,}ds=0,\qquad w^{i}\in W^{i},\;w^{j}\in W^{j}
\label{eq:starting_flux_constraint}%
\end{equation}
where $\mathbf{n}^{i}$ denotes the unit outward normal of $\Omega^{i}$.

\begin{assumption}
\label{ass:enough-constraints}Initial flux constraints
(\ref{eq:starting_flux_constraint}) are prescribed over all faces.
\end{assumption}

This set of initial constraints will be enriched by the adaptive method
described in Section~\ref{sec:adaptive}. Now, let us define $\widetilde{W}%
_{\Pi}\subset\widetilde{W}$ as the subspace of functions with values given by
the flux coarse degrees of freedom between adjacent substructures, and such
that they are Stokes harmonic, and let us also define $\widetilde{W}_{\Delta
}\subset\widetilde{W}$ as the subspace of all function such that
their flux coarse degrees of freedom vanish. The functions in$~\widetilde{W}%
_{\Pi}$ are uniquely determined by the values of their coarse degrees of
freedom, and%
\begin{equation}
\widetilde{W}=\widetilde{W}_{\Delta}\oplus\widetilde{W}_{\Pi}.
\label{eq:tilde-dec}%
\end{equation}

The next ingredient is the projection $E:\widetilde{W}\rightarrow\widehat{W}$
defined by taking a weighted average of corresponding degrees of freedom on
substructure interfaces, cf. Remark~\ref{rem:scaling}.

In implementation, we define$~\widetilde{W}$ using a matrix$~C_{U}$, which is
a block diagonal with blocks$~C_{U}^{i}$, $i=1,\dots,N$, and it is constructed
exactly as matrix~$C$ in~\cite[Section~2.3]{Mandel-2007-ASF},
\begin{equation}
\widetilde{W}=\left\{  w\in\left(  I-P\right)  W:C_{U}\left(  I-E\right)
w=0\right\}  . \label{eq:tilde-def}%
\end{equation}
The values $C_{U}v$ will be called local flux coarse degrees of freedom, and
the space$~\widetilde{W}$ consists of all functions such that their flux
coarse degrees of freedom on adjacent substructures have zero jumps. The
decomposition of the space $Q_{I}$ given by~(\ref{eq:Q_decomposition}) can be
also managed by\ constraints. We remark that this is somewhat non-standard
practice in substructuring, because the constraints are commonly related only
to the degrees of freedom at the interfaces. So, we define a space$~Q^{i}$,
for $i=1,\dots,N$, as
\begin{equation}
Q^{i}=\left\{  \left(  q\in Q\right)  |_{\Omega^{i}}:C_{Q}^{i}q=0\right\}  ,
\label{eq:Qi-def}%
\end{equation}
where the matrices~$C_{Q}^{i}$ are selected so
that~(\ref{eq:Q_decomposition-int}) is satisfied. In implementation,~$C_{Q}%
^{i}$ is a row vector with entries given by volumes of finite elements in
subdomain~$i$. Now we have all ingredients to recall the two-level BDDC
method~\cite[Algorithm~2]{Sousedik-2013-NBS}.

\begin{algorithm}
[BDDC method]\label{alg:two-level-nested}Find the solution $\left(
u,p\right)  \in\left(  U,Q\right)  $\ of problem (\ref{eq:variational-1}%
)--(\ref{eq:variational-2}) by computing:

\begin{enumerate}
\item the coarse component $u_{0}\in\widehat{W}$: solving $\left(
\widetilde{w}_{0},p_{0}\right)  \in\left(  \widetilde{W}_{\Pi},Q_{0}\right)  $
from
\begin{align}
a\left(  \widetilde{w}_{0},\widetilde{v}_{\Pi}\right)  +b\left(
\widetilde{v}_{\Pi},p_{0}\right)   &  =0,\qquad\forall\widetilde{v}_{\Pi}%
\in\widetilde{W}_{\Pi},\label{eq:two-level-nested_coarse-1}\\
b\left(  \widetilde{w}_{0},q_{0}\right)   &  =\left\langle f,q_{0}%
\right\rangle ,\qquad\forall q_{0}\in Q_{0},
\label{eq:two-level-nested_coarse-2}%
\end{align}
dropping $p_{0}$, and applying the projection
\[
u_{0}=E\widetilde{w}_{0}.
\]

\item the substructure components $\left(  u_{I},p_{I}\right)  \in\left(
U_{I},Q_{I}\right)  $ from
\begin{align*}
a\left(  u_{I},v_{I}\right)  +b\left(  v_{I},p_{I}\right)   &  =-a\left(
u_{0},v_{I}\right)  ,\qquad\forall v_{I}\in U_{I},\\
b\left(  u_{I},q_{I}\right)   &  =\left\langle f,q_{I}\right\rangle -b\left(
u_{0},q_{I}\right)  ,\qquad\forall q_{I}\in Q_{I},
\end{align*}
\qquad\qquad

dropping $p_{I}$, and adding the solutions as
\begin{equation}
u^{\ast}=u_{0}+u_{I}. \label{eq:two-level_u^star}%
\end{equation}

\item the correction and the pressure $\left(  u_{\text{corr}},p\right)
\in\left(  U,Q\right)  $ from%
\begin{align}
a\left(  u_{\text{corr}},v\right)  +b\left(  v,p\right)   &  =-a\left(
u^{\ast},v\right)  ,\qquad\forall v\in U,\label{eq:two-level_corr-1}\\
b\left(  u_{\text{corr}},q\right)   &  =0,\qquad\forall q\in Q.
\label{eq:two-level_corr-2}%
\end{align}
Specifically, use the CG method with the BDDC\ preconditioner defined in
Algorithm~\ref{alg:two-level}, using the same setup of the coarse problem as
in (\ref{eq:two-level-nested_coarse-1})--(\ref{eq:two-level-nested_coarse-2}).

Finally, the flux variables are obtained as
\[
u=u^{\ast}+u_{\text{corr}}.
\]

\end{enumerate}
\end{algorithm}

\begin{remark}
\label{rem:scaling}The difference between problems (\ref{eq:variational-1}%
)--(\ref{eq:variational-2}) and (\ref{eq:two-level_corr-1}%
)--(\ref{eq:two-level_corr-2}) is that the latter problem has a vanishing
second component, and therefore the correction $u_{\text{corr}}$ is
divergence-free by (\ref{eq:two-level_corr-2}). Also, we note that the initial
flux constraints constructed according to~(\ref{eq:starting_flux_constraint})
do not allow scaling weights in the scaling operator$~E$ to vary along the
interface in order for $u^{\ast}$ to satisfy
\[
b\left(  u^{\ast},q\right)  =\left\langle f,q\right\rangle ,\qquad\forall q\in
Q.
\]
Therefore, in our numerical experiments, we use the multiplicity scaling
unless the coefficient jumps are aligned with subdomain interfaces, see
also~\cite[Remark~2]{Sousedik-2013-NBS}.
\end{remark}

The application of the BDDC\ preconditioner for the computation of $u_{corr}$
using two- resp. three-level method was studied by
Tu~\cite{Tu-2005-BAM,Tu-2011-TBA}. In~\cite{Sousedik-2013-NBS}, we applied
Algorithm~\ref{alg:two-level-nested} recursively. Here, we will introduce a
specific construction of the space$~\widetilde{W}_{\Pi}$ but before doing so,
let us discuss Step~3 of Algorithm~\ref{alg:two-level-nested} in more detail.

The first step in substructuring is typically reduction of the problem to
interfaces. In particular, problem~(\ref{eq:two-level_corr-1}%
)--(\ref{eq:two-level_corr-2}) is reduced to finding $\left(  \widehat{w}%
,p_{0}\right)  \in\left(  \widehat{W},Q_{0}\right)  $ such that%
\begin{align}
a\left(  \widehat{w},\widehat{v}\right)  +b\left(  \widehat{v},p_{0}\right)
&  =\left\langle f^{\ast},\widehat{v}\right\rangle ,\qquad\forall
\widehat{v}\in\widehat{W},\label{eq:corr-reduced-1}\\
b\left(  \widehat{u},q_{0}\right)   &  =0,\qquad\forall q_{0}\in Q_{0},
\label{eq:corr-reduced-2}%
\end{align}
where $f^{\ast}\in\widehat{W}^{\prime}$ is the reduced right-hand side. In
implementation, the interiors are eliminated by the static condensation,
problem~(\ref{eq:corr-reduced-1})--(\ref{eq:corr-reduced-2}) is solved
iteratively, and the interiors $\left(  u_{I},p_{I}\right)  \in\left(
U_{I},Q_{I}\right)  $ are recovered in the post-correction. The key
observation is, cf.~\cite[Section 9.4.2]{Toselli-2005-DDM}, that if we define
a \emph{balanced} subspace%
\[
\widehat{W}_{B}=\left\{  \widehat{w}\in\widehat{W}:b\left(  \widehat{w}%
,q_{0}\right)  =0,\quad\forall q_{0}\in Q_{0}\right\}  ,
\]
problem~(\ref{eq:corr-reduced-1})--(\ref{eq:corr-reduced-2}) becomes
equivalent to the positive definite problem
\[
\widehat{u}\in\widehat{W}_{B}:\quad a\left(  \widehat{u},\widehat{v}\right)
=\left\langle f^{\ast},v\right\rangle ,\quad\forall\widehat{v}\in
\widehat{W}_{B}.
\]
This observation justifies use of the CG method\ preconditioned by the BDDC
provided that an initial guess is balanced, e.g., zero, and the outputs of the
preconditioner are also balanced.
It also implies that the iterates are effectively performed with the flux
unknowns, and the pressure components$~p_{0}$ are resolved in the coarse
correction of the preconditioner. The precise formulation of the two-level
BDDC\ preconditioner for saddle-point problems follows. It is the reduced
variant of~\cite[Algorithm~3]{Sousedik-2013-NBS}.

\begin{algorithm}
[BDDC preconditioner]\label{alg:two-level} Define the preconditioner $\left(
r_{B},0\right)  \in\left(  \widehat{W}^{\prime},Q_{0}^{\prime}\right)
\longmapsto\left(  \widehat{w},p_{0}\right)  \in\left(  \widehat{W}%
,Q_{0}\right)  $\ by computing:

\begin{enumerate}
\item the coarse correction $\left(  w_{\Pi},p_{0}\right)  \in\left(
\widetilde{W}_{\Pi},Q_{0}\right)  $ from
\begin{align}
a\left(  w_{\Pi},z_{\Pi}\right)  +b\left(  z_{\Pi},p_{0}\right)   &
=\left\langle r_{B},Ez_{\Pi}\right\rangle ,\qquad\forall z_{\Pi}%
\in\widetilde{W}_{\Pi},\label{eq:two-level_coarse-1}\\
b\left(  w_{\Pi},q_{0}\right)   &  =0,\qquad\forall q_{0}\in Q_{0}.
\label{eq:two-level_coarse-2}%
\end{align}

\item the substructure correction $w_{\Delta}\in\widetilde{W}_{\Delta}$ from%
\begin{align*}
a\left(  w_{\Delta},z_{\Delta}\right)  +b\left(  z_{\Delta},p_{I\Delta
}\right)   &  =\left\langle r_{B},Ez_{\Delta}\right\rangle ,\qquad\forall
z_{\Delta}\in\widetilde{W}_{\Delta},\\
b\left(  w_{\Delta},q_{I}\right)   &  =0,\qquad\forall q_{I}\in Q_{I}.
\end{align*}

\item the sum and average of the two corrections%
\begin{equation}
\widehat{w}=E\left(  w_{\Pi}+w_{\Delta}\right)  . \label{eq:two-level_w}%
\end{equation}

\end{enumerate}
\end{algorithm}

In order to state the condition number bound, we also need to introduce a
larger space of \emph{balanced} functions $\widetilde{W}_{B}$ such that
$\widehat{W}_{B}\subset\widetilde{W}_{B}$\ defined as
\[
\widetilde{W}_{B}=\left\{  w\in\widetilde{W}:b\left(  v,q_{0}\right)
=0,\quad\forall q_{0}\in Q_{0}\right\}  .
\]
The space$~\widetilde{W}_{\Pi}$ is also balanced, i.e., $\widetilde{W}_{\Pi
}\subset\widetilde{W}_{B}$ by (\ref{eq:two-level_coarse-2}). Then also the
output of the preconditioner~(\ref{eq:two-level_w}) satisfies $\widehat{w}%
\in\widehat{W}_{B}$, and we refer to~\cite[Lemma~3]{Sousedik-2013-NBS} for the proof.

Finally, we formulate the condition number bound. If we note that $E$ is a
projection, it is the same as~\cite[Theorem~4]{Sousedik-2013-NBS}
or~\cite[Theorem~6.1]{Tu-2005-BAM}, cf. also~\cite[Theorem$~$3]%
{Mandel-2007-ASF}.

\begin{theorem}
\label{thm:two-level-bound}The condition number~$\kappa$ of the BDDC
preconditioner from Algorithm~\ref{alg:two-level} satisfies
\begin{equation}
\kappa\leq\omega=\max\left\{  {\sup_{w\in\widetilde{W}_{B}}\frac{\left\Vert
\left(  I-E\right)  w\right\Vert _{a}^{2}}{\left\Vert w\right\Vert _{a}^{2}%
},1}\right\}  \leq C\left(  1+\log\frac{H}{h}\right)  ^{2}\text{.}
\label{eq:two-level-bound}%
\end{equation}

\end{theorem}

The bound$~\omega$ in~(\ref{eq:two-level-bound}) inspires the adaptive
selection of the flux constraints.

\section{Adaptive selection of the flux constraints}

\label{sec:adaptive}The basic idea is same as in our previous work on adaptive
BDDC for elliptic
problems~{\cite{Mandel-2007-ASF,Mandel-2012-ABT,Sousedik-2013-AMB}}. The
bound~$\omega$ in~(\ref{eq:two-level-bound}) is equal to the maximal
eigenvalue $\lambda_{\max}$\ of the generalized eigenvalue problem%
\begin{equation}
w\in\widetilde{W}_{B}:\quad a\left(  \left(  I-E\right)  w,\left(  I-E\right)
z\right)  =\lambda\,a\left(  w,z\right)  ,\quad\forall z\in\widetilde{W}_{B}.
\label{eq:eigenvalue-problem}%
\end{equation}
From the Courant-Fisher-Weyl minimax principle cf., e.g.,~\cite[Theorem
5.2]{Demmel-1997-ANL}, the bound$~\omega$ can be decreased by adding
constraints in the definition of the space$~\widetilde{W}_{B}$ as:

\begin{lemma}
[\cite{Mandel-2012-ABT,Sousedik-2013-AMB}]The generalized eigenvalue
problem~(\ref{eq:eigenvalue-problem}) has eigenvalues $\lambda_{1}\geq
\lambda_{2}\geq\dots\lambda_{n}\geq0$. Denote the corresponding eigenvectors
by$~w_{\ell}$. Then, for any $k=1,\dots,n-1$, and any linear functionals
$L_{\ell}$, $\ell=1,\dots,k$,%
\[
\max\left\{  \frac{\left\Vert \left(  I-E\right)  w\right\Vert _{a}^{2}%
}{\left\Vert w\right\Vert _{a}^{2}}:w\in\widetilde{W}_{B},\;L_{\ell}\left(
w\right)  =0\;\;\forall\ell=1,\dots,k\right\}  \geq\lambda_{k+1},
\]
with equality if%
\[
L_{\ell}\left(  w\right)  =a\left(  \left(  I-E\right)  w_{\ell},\left(
I-E\right)  w\right)  .
\]

\end{lemma}

Because solving the global eigenvalue problem\ (\ref{eq:eigenvalue-problem})
is computationally expensive, we replace it by a collection of much smaller
problems defined for all pairs of adjacent substructures, where a pair of
substructures is adjacent if they share a face. All quantities associated with
a pair of adjacent substructures $\Omega^{i}$ and $\Omega^{j}$ will be denoted
by a superscript$~^{ij}$. In particular, we define $W^{ij}=W^{i}\times W^{j}$,
and the local space~$\widetilde{W}^{ij}$ of Stokes harmonic functions that
satisfy the initial constraints at the face$~\Gamma^{ij}$ by%
\begin{equation}
\widetilde{W}_{B}^{ij}=\left\{  w\in\left(  I-P^{ij}\right)  W^{ij}:C_{U}%
^{ij}\left(  I-E^{ij}\right)  w=0\right\}  . \label{eq:def-Wtilde-dual}%
\end{equation}
We note that the space$~\widetilde{W}_{B}^{ij}$ is balanced, which is an
implication of Assumption~\ref{ass:enough-constraints}.

In these settings~(\ref{eq:eigenvalue-problem}) becomes a \emph{local} problem
to find $w\in\widetilde{W}_{B}^{ij}$ such that
\begin{equation}
a^{ij}\left(  \left(  I-E^{ij}\right)  w,\left(  I-E^{ij}\right)  z\right)
=\lambda\,a^{ij}\left(  w,z\right)  ,\quad\forall z\in\widetilde{W}_{B}^{ij}.
\label{eq:eigenvalue-problem-local}%
\end{equation}
The bilinear form $a^{ij}$\ is associated on$~\widetilde{W}_{B}^{ij}$ with the
Schur complement$~S^{ij}$ defined with respect to the interfaces$~\Gamma^{i}%
$,$~\Gamma^{j}$, and is positive-definite, cf.~\cite[Lemma$~$3.1]{Tu-2005-BAM}.

Now we can proceed in the same way as in~\cite{Sousedik-2013-AMB}. Let us
denote by$~\mathcal{C}$ the matrix corresponding to $C_{U}^{ij}\left(
I-E^{ij}\right)  $. The orthogonal projection onto null$~\mathcal{C}$ is given
by%
\[
\Pi=I-\mathcal{C}^{T}\left(  \mathcal{CC}^{T}\right)  ^{-1}\mathcal{C},
\]
and we implement the local generalized eigenvalue
problems~(\ref{eq:eigenvalue-problem-local}) as%
\begin{equation}
\Pi\left(  I-E^{ij}\right)  ^{T}S^{ij}\left(  I-E^{ij}\right)  \Pi
w=\lambda\,\Pi S^{ij}\Pi w, \label{eq:eigenvalue-problem-local-matrix}%
\end{equation}
which can be either solved using a dense eigenvalue
solver~\cite{Mandel-2007-ASF} or eventually, since
\[
\text{null}\left[  \Pi S^{ij}\Pi\right]  \subset\text{null}\left[  \Pi\left(
I-E^{ij}\right)  ^{T}S^{ij}\left(  I-E^{ij}\right)  \Pi\right]  ,
\]
a subspace iterations such as the LOBPCG method~\cite{Knyazev-2001-TOP}, which
runs effectively in the factorspace, could be also used.
From~(\ref{eq:eigenvalue-problem-local-matrix}), we wish the constraints to
satisfy
\[
L_{\ell}\left(  w\right)  =w_{\ell}^{T}\Pi\left(  I-E^{ij}\right)  ^{T}%
S^{ij}\left(  I-E^{ij}\right)  \Pi w=0.
\]
That is, we would add into the matrix $C_{U}^{ij}$ the rows
\begin{equation}
c_{\ell}^{ij}=w_{\ell}^{T}\Pi\left(  I-E^{ij}\right)  ^{T}S^{ij}\left(
I-E^{ij}\right)  \Pi, \label{eq:constraints}%
\end{equation}
but because by~\cite[Proposition 1]{Sousedik-2013-AMB} each row can be split
as $c_{\ell}^{ij}=\left[
\begin{array}
[c]{cc}%
c_{\ell}^{i} & -c_{\ell}^{i}%
\end{array}
\right]  $ and either half of$~c_{\ell}^{ij}$\ is used to augment the
matrices$~C_{U}^{i}$ and$~C_{U}^{j}$, see~(\ref{eq:C^i}). We note that, due to
the discretization using RT0 elements, the added rows are readily available in
the form used in substructuring. The adaptive BDDC algorithm follows.

\begin{algorithm}
[Adaptive BDDC]\label{alg:adaptive}Find the smallest $k$ for every two
adjacent substructures to guarantee that $\lambda_{k+1}\leq\tau$, where $\tau$
is a given tolerance threshold (the target condition number), and add the
constraints~(\ref{eq:constraints}) to the definition of $\widetilde{W}$.
\end{algorithm}

After the adaptive constraints are added, we define the \emph{heuristic
condition number indicator} as the largest eigenvalue $\omega^{ij}$ of all
local eigenvalue problems~(\ref{eq:eigenvalue-problem-local}), that is
\begin{equation}
\widetilde{\omega}=\max\left\{  \omega^{ij}:\Omega^{i}\text{ and }\Omega
^{j}\text{ are adjacent}\right\}  . \label{eq:indicator}%
\end{equation}

\begin{remark}
It has been shown in~\cite[Theorem~4.3]{Zampini-2017-MBD}, see
also~\cite[Theorem~3.10]{Pechstein-2017-UFA} and~\cite[Theorem~3.3]%
{Oh-2018-BAD}, that the condition number$~\kappa$ of the adaptive BDDC
operator satisfies
\[
\kappa\leq\widetilde{\omega}\,N_{F}^{2},
\]
where $N_{F}$ is the maximum number of faces of any subdomain. We note that
this bound is pessimistic due to the factor$~N_{F}^{2}$, and in fact we
observed $\kappa\approx\widetilde{\omega}$ in all experiments.
\end{remark}

\section{Implementation remarks}

\label{sec:implementation}First, we describe a rescaling used to preserve
numerical stability of the method with highly heterogeneous permeability
coefficients. The variational problem~(\ref{eq:variational-1}%
)--(\ref{eq:variational-2}) can be written in the matrix form as%
\begin{equation}
\left[
\begin{array}
[c]{cc}%
A & B^{T}\\
B & 0
\end{array}
\right]  \left[
\begin{array}
[c]{c}%
u\\
p
\end{array}
\right]  =\left[
\begin{array}
[c]{c}%
0\\
f
\end{array}
\right]  . \label{eq:system}%
\end{equation}
Assuming that the mesh size$~h\approx1$, the entries in$~A$ are$~O\left(
k^{-1}\right)  $ and the entries in$~B$ are$~O\left(  1\right)  $. In
particular, in the case of the SPE~10 data set we get $k^{-1}\approx
10^{6}-10^{12}$, and we found that some of the subdomain matrices and the
matrix of the coarse problem may appear numerically singular. Due to the
discontinuous approximation of the pressure, $B$ is a block-diagonal
rectangular matrix. Each block corresponds to a particular subdomain$,$ and it
can be rescaled, e.g., by an average of the diagonal entries of $A$
corresponding to the degrees of freedom in this subdomain. Collecting this
scaling coefficients in a diagonal matrix$~D$, we replace~(\ref{eq:system})
by
\begin{equation}
\left[
\begin{array}
[c]{cc}%
A & B^{T}D\\
DB & 0
\end{array}
\right]  \left[
\begin{array}
[c]{c}%
u\\
\overline{p}%
\end{array}
\right]  =\left[
\begin{array}
[c]{c}%
0\\
Df
\end{array}
\right]  , \label{eq:system-scaled}%
\end{equation}
and the pressure is recovered at the end of computations as $p=D\overline{p}$.

\subsection{Coarse degrees of freedom}

\label{sec:coarse-dofs}The selection of the flux coarse degrees of freedom or,
equivalently, flux constraints entails construction of the matrix~$C_{U}$ in
the definition of the space~$\widetilde{W}$ by~(\ref{eq:tilde-def}).
Similarly, the selection of the pressure constraints, which facilitate the
decomposition~(\ref{eq:Q_decomposition}), entails construction of the
matrices~$C_{Q}^{i}$, $i=1,\dots,N$, in the definition of the spaces$~Q^{i}$
by~(\ref{eq:Qi-def}). Following the standard practice in substructuring, in
implementation we work with global and local degrees of freedom and the
corresponding spaces, and vectors from these spaces are related by a
restriction operator (a~zero-one matrix). Therefore, the matrix$~C_{U}$ is
constructed as a block-diagonal matrix using blocks$~C_{U}^{i}$ that select
local flux coarse degrees of freedom from all degrees of freedom of
substructure~$i$, see~\cite[Section~2.3]{Mandel-2007-ASF} for details. In the
mixed finite element settings, each local coarse degrees of freedom selection
matrix is constructed simply by augmenting the matrix$~C_{U}^{i}$ by a
row$~C_{Q}^{i}$ as
\begin{equation}
\left[
\begin{array}
[c]{cc}%
C_{U}^{i} & \\
& C_{Q}^{i}%
\end{array}
\right]  ,\qquad i=1,\dots,N, \label{eq:C^i}%
\end{equation}
and the matrices$~C_{U}^{i}$ may be further augmented by the adaptive
algorithm, see~(\ref{eq:constraints}).

\subsection{Solution of the local generalized eigenvalue problems}

The choice of an eigensolver for the eigenvalue
problems~(\ref{eq:eigenvalue-problem-local-matrix}) is a delicate one. In
general, the decision whether to use a dense or a sparse\ eigensolver depends
on the type of the eigenvalue problem, size of the substructures, dimension of
the problem, availability of a preconditioner for a sparse solver, and
conditioning and numerical sensitivity of the underlying problem. All these
factors will clearly affect the overall computational cost and performance of
the method. We note that the
formulation~(\ref{eq:eigenvalue-problem-local-matrix}) allows to use a
matrix-free iterative method such as the LOBPCG~\cite{Knyazev-2001-TOP} in the
same way as for elliptic problems, including that it can be further
preconditioned by a local version of the BDDC as suggested in~\cite[Section~5]%
{Sousedik-2013-AMB}, see also~\cite{Klawonn-2018-CLL}. However, we found
that\ dense eigenvalue solvers are more suitable for the SPE~10 dataset due to
their robustness, and we used \textsc{Matlab} function \texttt{eig} in the
numerical experiments.

\subsection{Computational cost}

Clearly, the two most computationally expensive parts of the method are the
setup of the constraints by solving the set of the local eigenvalue problems,
and the factorization of the coarse problem. There are many eigenvalue
problems to be solved, but they are small and can be solved in parallel---this
feature is similar to the setup of multiscale finite element
methods~\cite{Efendiev-2009-MFE}. Assuming that these can be solved
efficiently, the bottleneck in computations is the factorization of the coarse
problem. Specifically, it is crucial for the application of the method to
appropriately balance the effort in the preconditioner and the global linear
solver through a judicious choice of$~\tau$. This could be, for example,
achieved as follows: one can partition the domain into subdomains balancing
the sizes of subdomains and assuming certain size of the coarse problem (and
ideally also taking into account the coefficient jumps and minimizing the size
of interfaces), solve the set of local eigenvalue problems, and based on the
eigenvalues determine the number of additional adaptive constraints (and hence
the value of$~\widetilde{\omega}$) which minimize the work needed to factor
the coarse problem and the work needed by preconditioned conjugate gradients,
including the coarse problem back-substitutions, needed to reduce the error to
desired accuracy based on the well-known error reduction formula of conjugate
gradients see, e.g.,~\cite[Theorem~10.2.6]{Golub-1996-MAC}.

\section{Numerical experiments}

\label{sec:numerical}
\begin{figure}[ptbh]
\begin{center}%
\begin{tabular}
[c]{cc}%
\includegraphics[width=5.6cm]{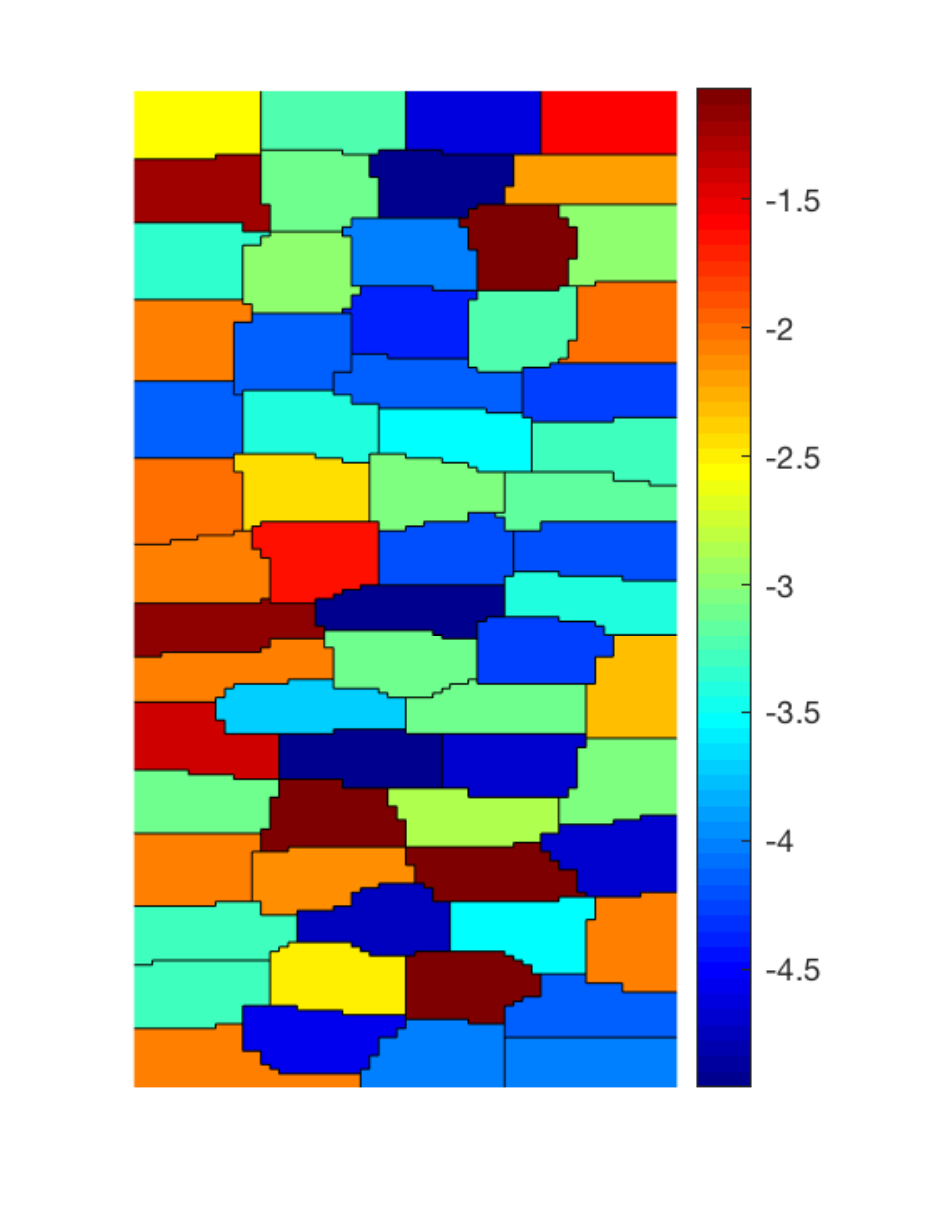} &
\raisebox{0.25\height}{\includegraphics[width=6.5cm]{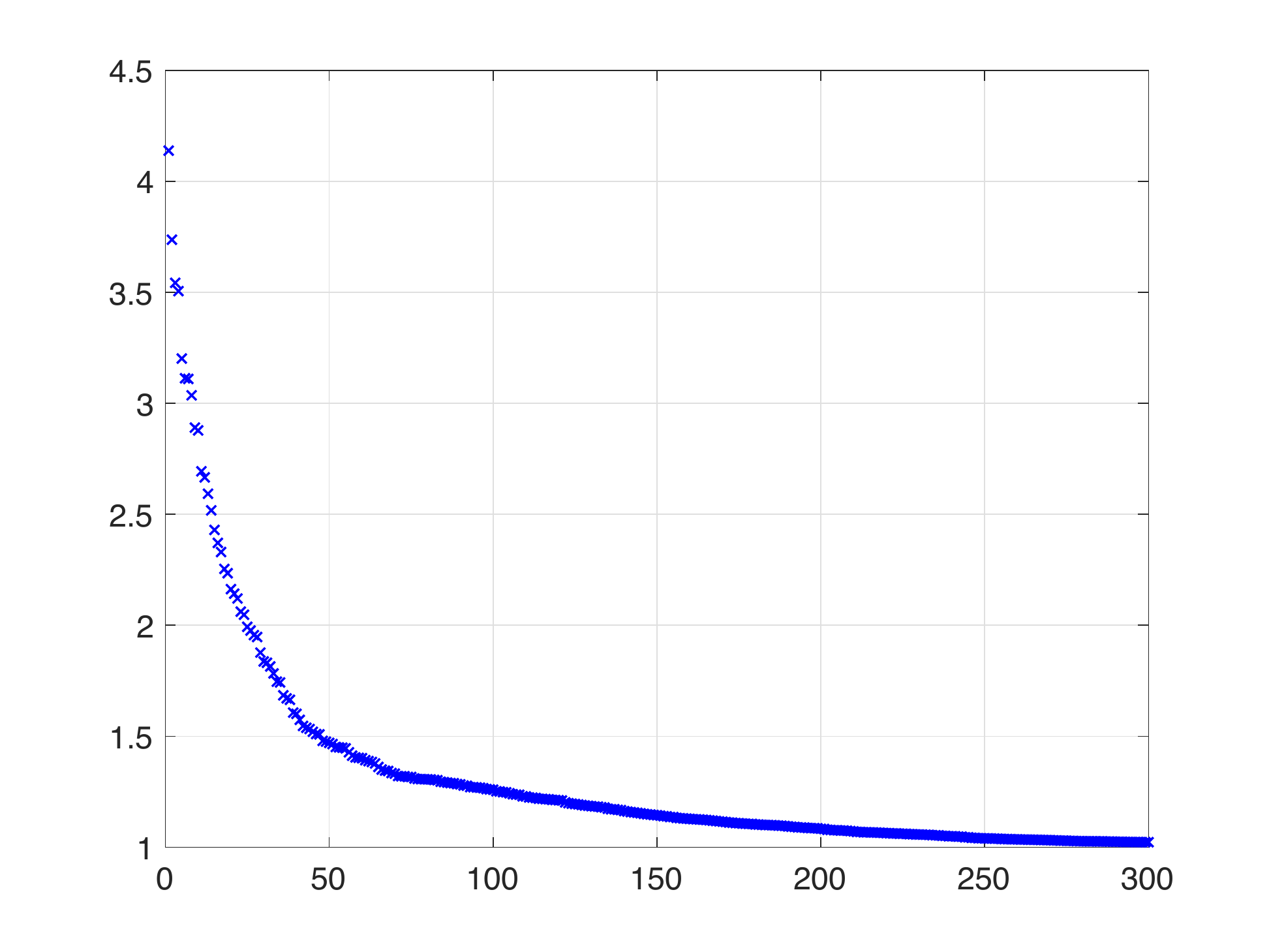}}
\end{tabular}
\end{center}
\caption{Substructuring and the base $10$ logarithm of the permeability $k$
for the problem with jumps aligned with the substructure interfaces (left
panel) and the largest $300$ eigenvalues of the BDDC\ preconditioned operator
for this problem (right panel). }%
\label{fig:jumps}%
\end{figure}

\begin{figure}[ptbh]
\begin{center}
\includegraphics[width=6.3cm]{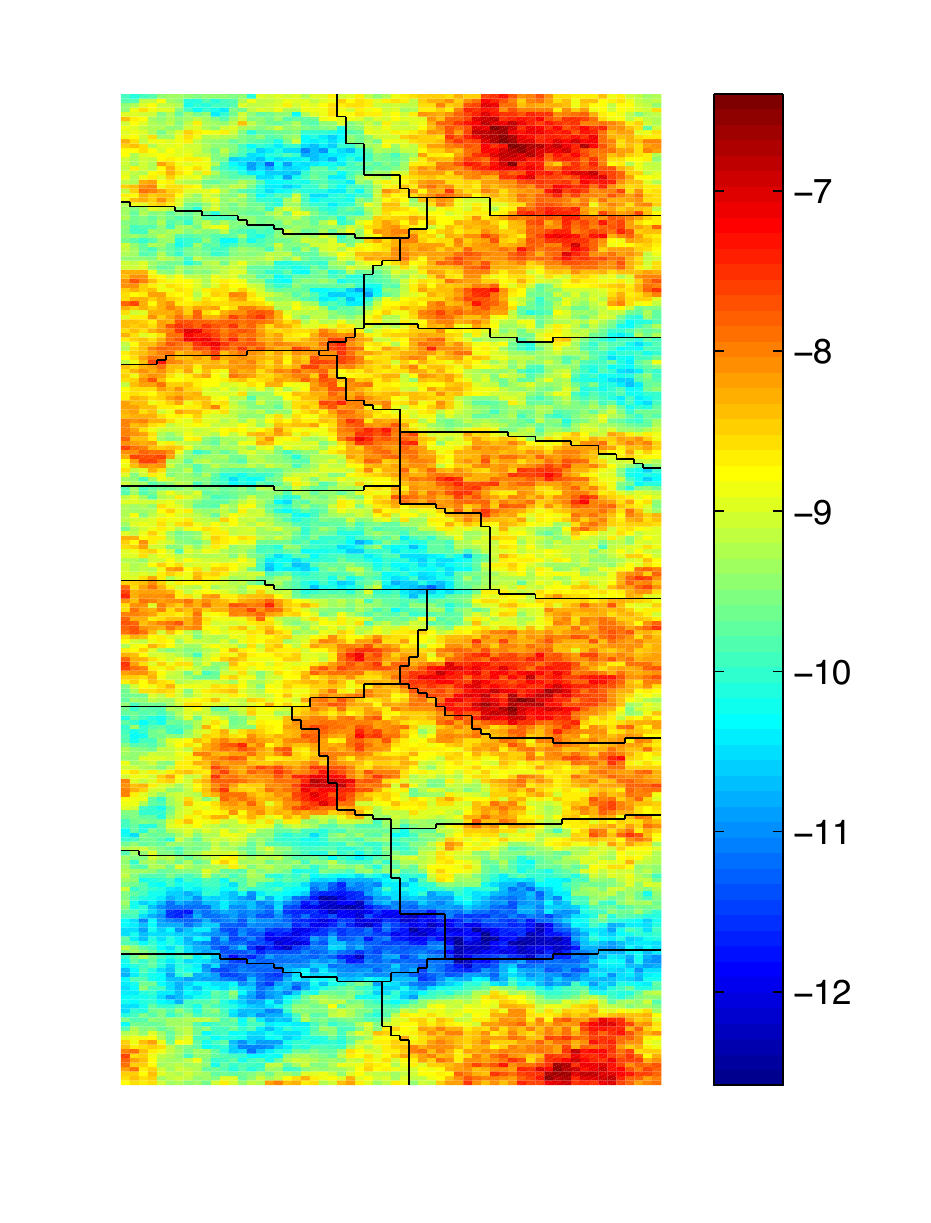}
\hspace{1mm}
\includegraphics[width=6.3cm]{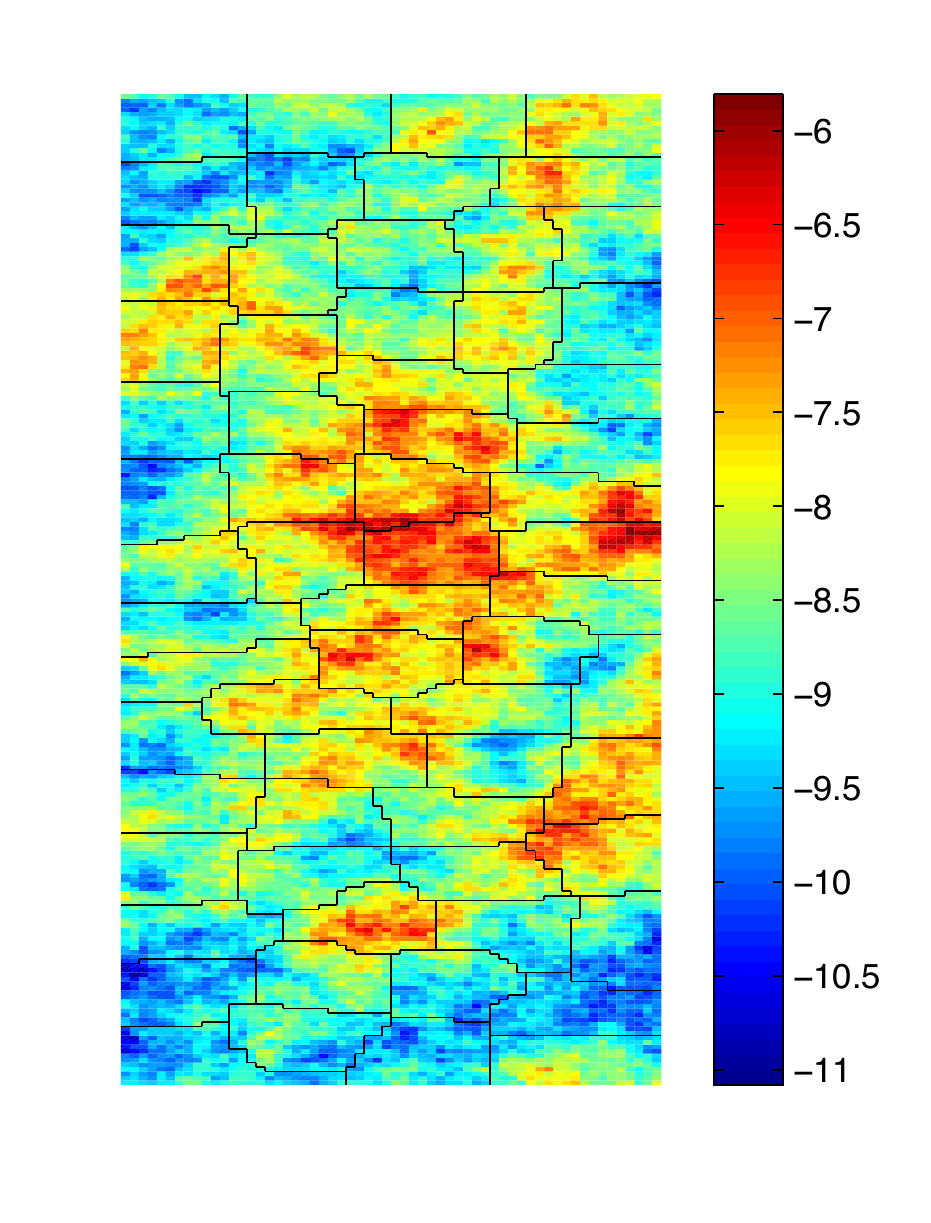}
\vspace{-10pt}
\end{center}
\caption{Substructuring and the base~$10$ logarithm of the permeability~$k$ in
layer~1 (left panel) and layer~20 (right panel) of the SPE~10 problem. Left
panel also illustrates irregular partitioning~A, and the right panel
illustrates irregular partitioning~B.}%
\label{fig:layers_1_20}%
\end{figure}

\begin{figure}[ptbh]
\begin{center}
\includegraphics[width=6.3cm]{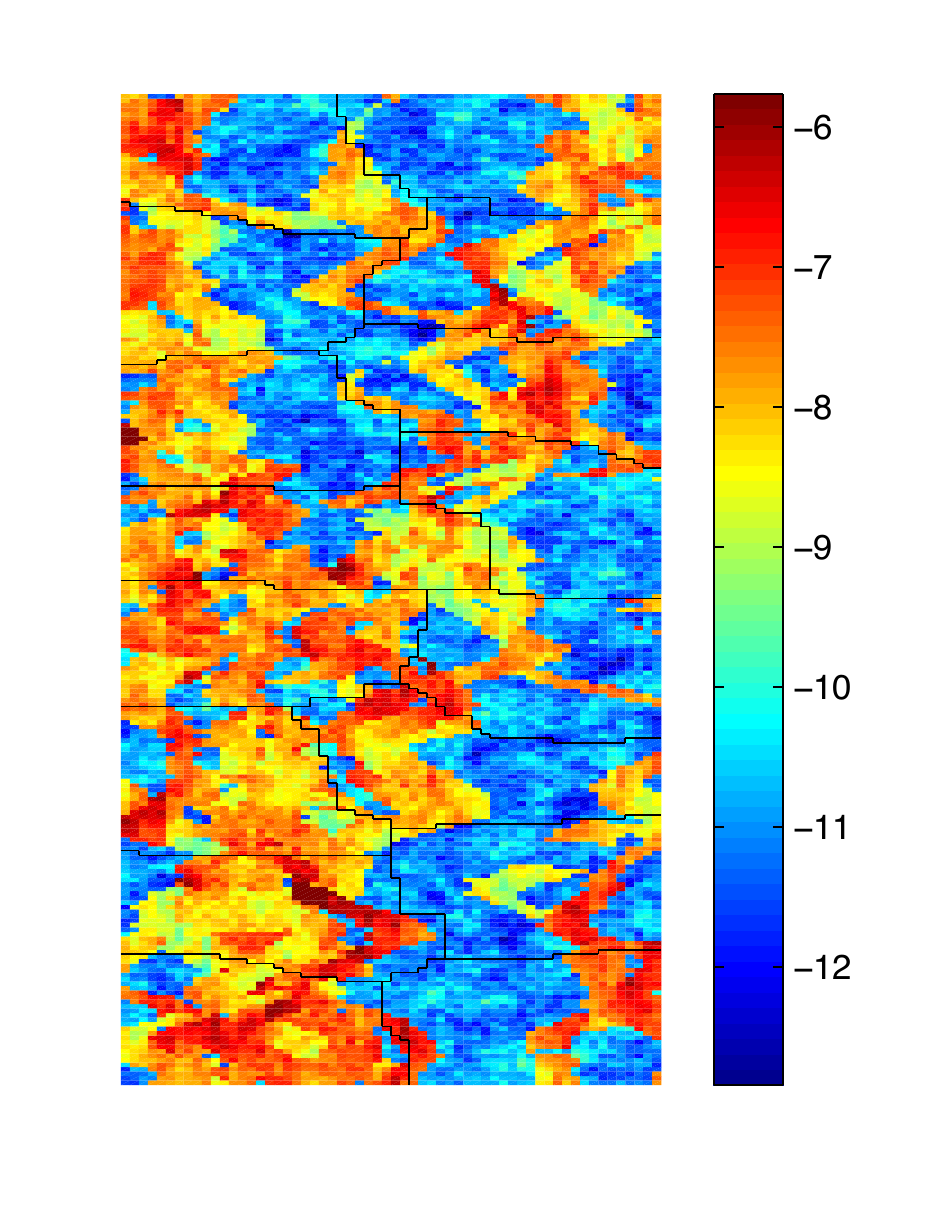}
\hspace{1mm}
\includegraphics[width=6.3cm]{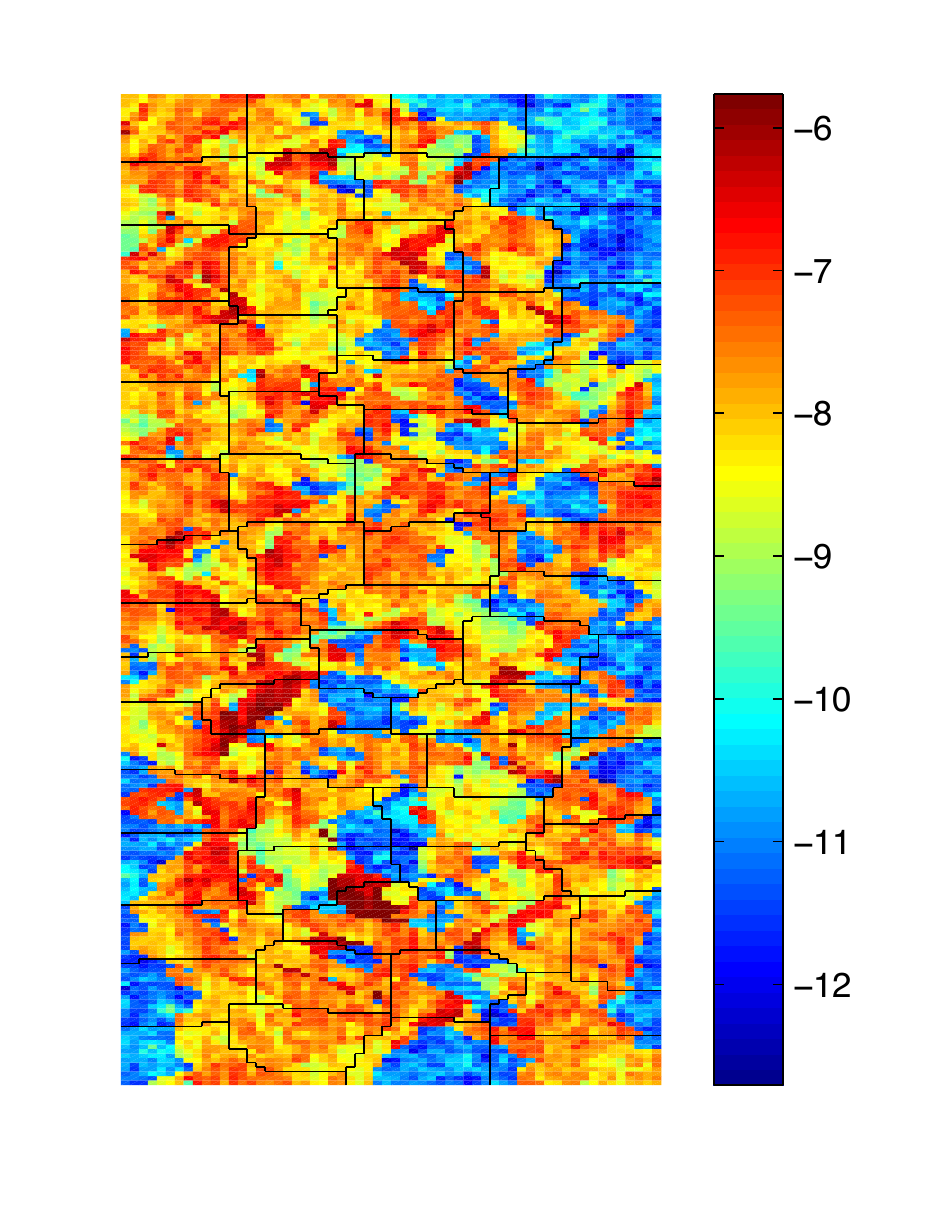}
\vspace{-10pt}
\end{center}
\caption{Substructuring and the base~$10$ logarithm of the permeability~$k$ in
layer~60 (left panel) and layer~85 (right panel) of the SPE~10 problem. Left
panel illustrates irregular partitioning A, and the right panel illustrates
irregular partitioning B.}%
\label{fig:layers_60_85}%
\end{figure}

\begin{figure}[ptbh]
\begin{center}
{\small permeability in layers~1--30:} \newline%
\includegraphics[width=6.4cm]{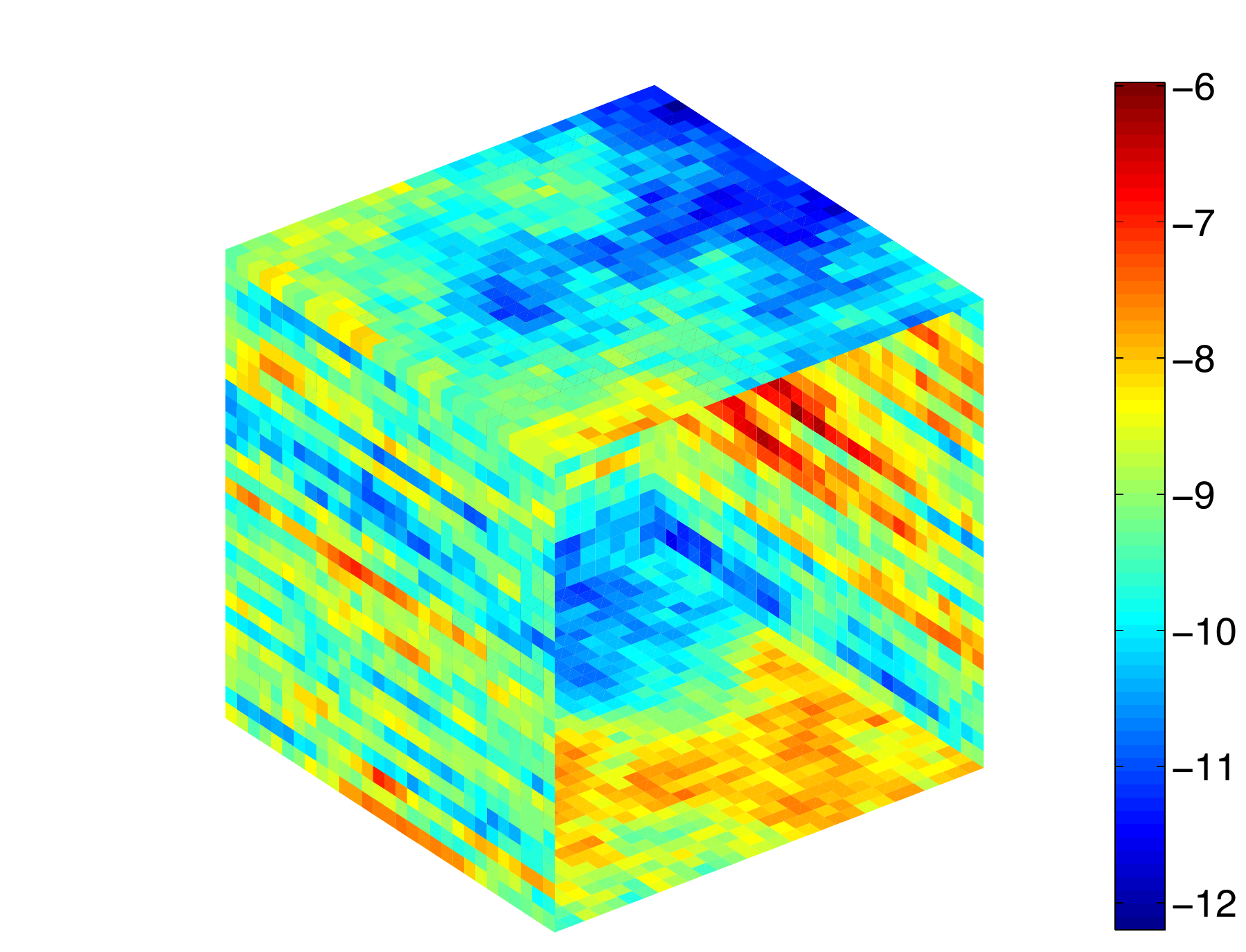}
\includegraphics[width=6.4cm]{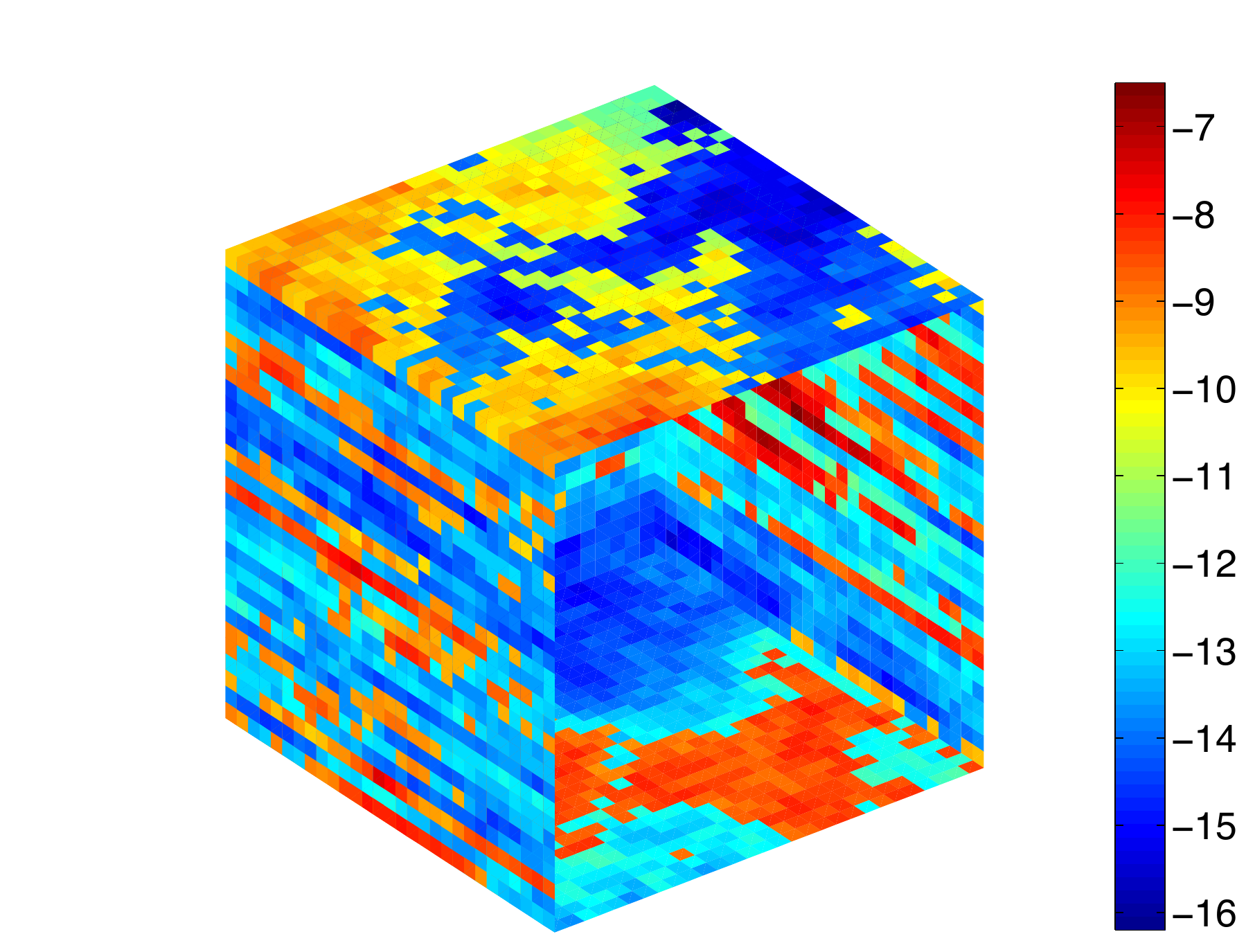}
\newline{\small permeability in layers~56--85:} \newline%
\includegraphics[width=6.4cm]{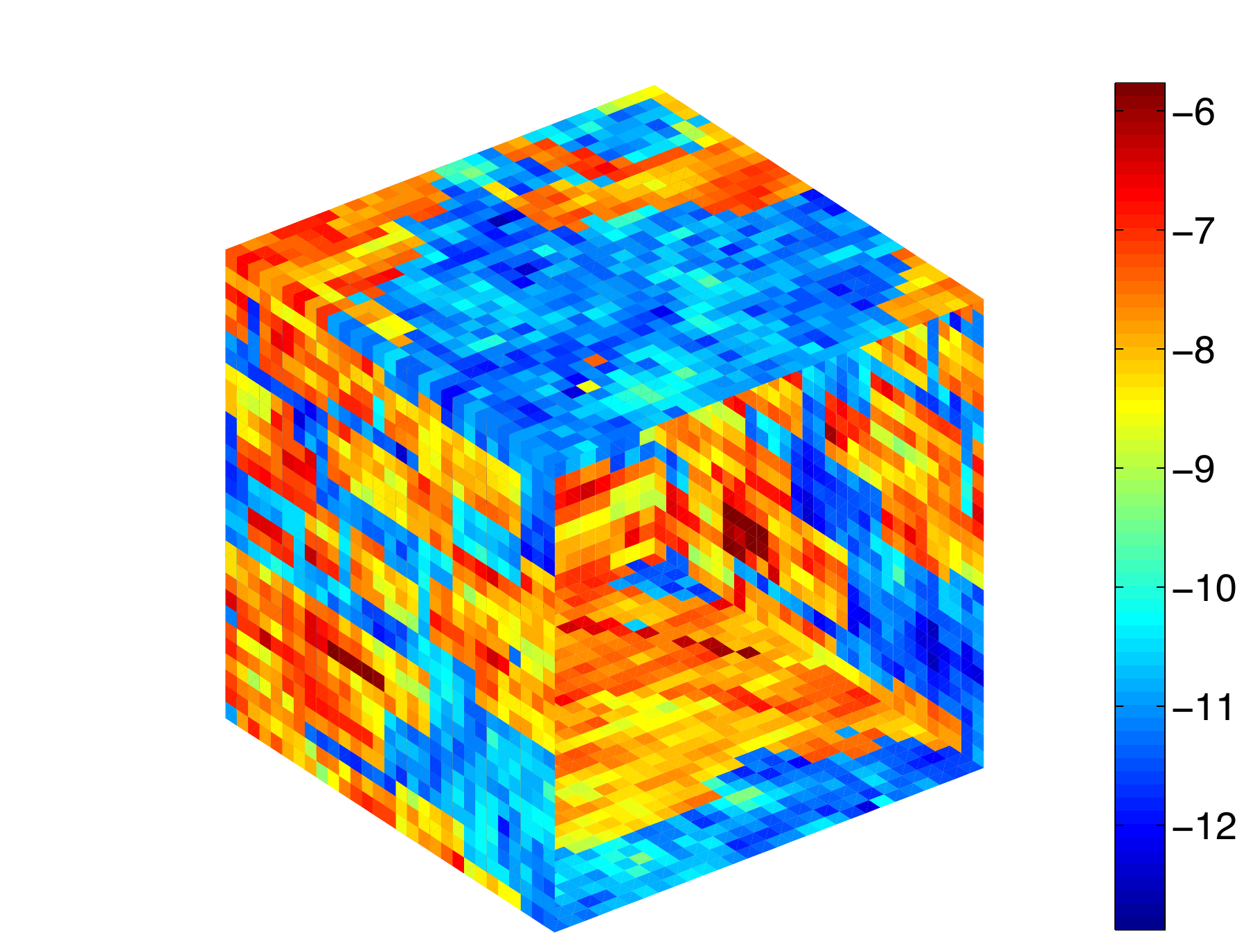}
\includegraphics[width=6.4cm]{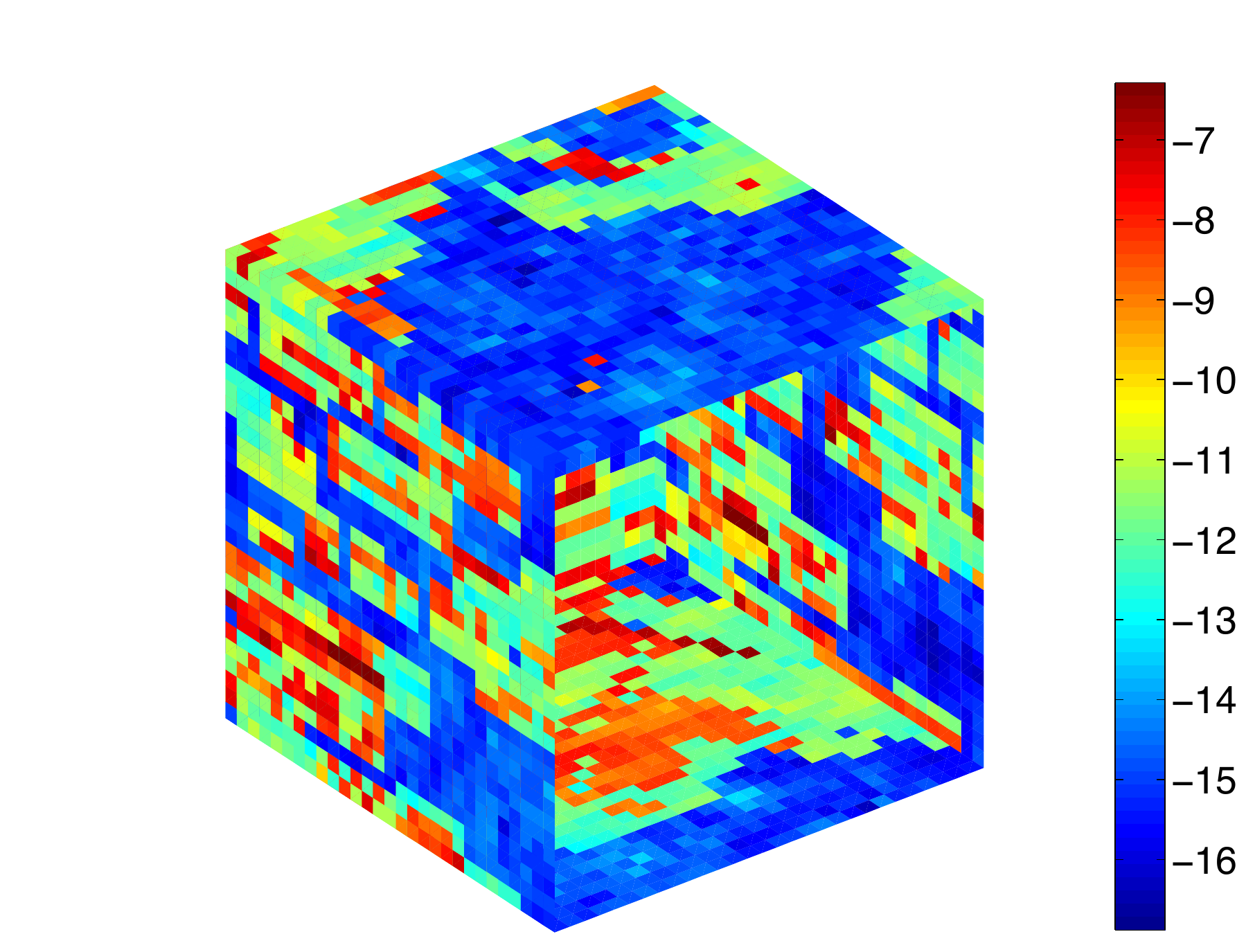}
\end{center}
\caption{Base $10$ logarithm of the permeability~$k$ in $x$ and $y$ directions
(left), and in $z$ direction (right) in two cutouts of the SPE~10 problem
consisting of $30\times30\times30$ elements.}%
\label{fig:mixed_RT0_3D}%
\end{figure}

\begin{figure}[ptbh]
\begin{center}
\includegraphics[width=6.4cm]{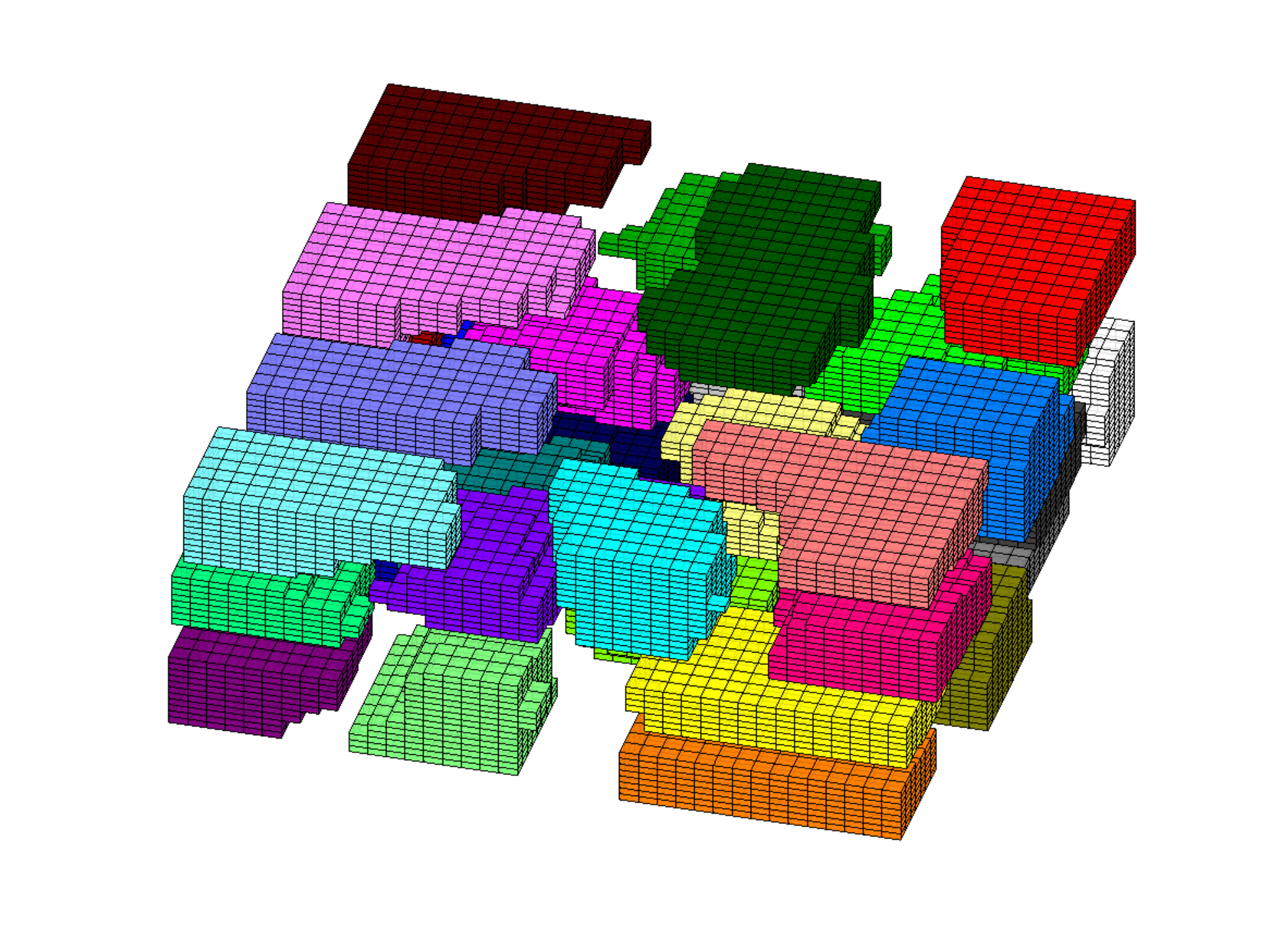}
\includegraphics[width=6.4cm]{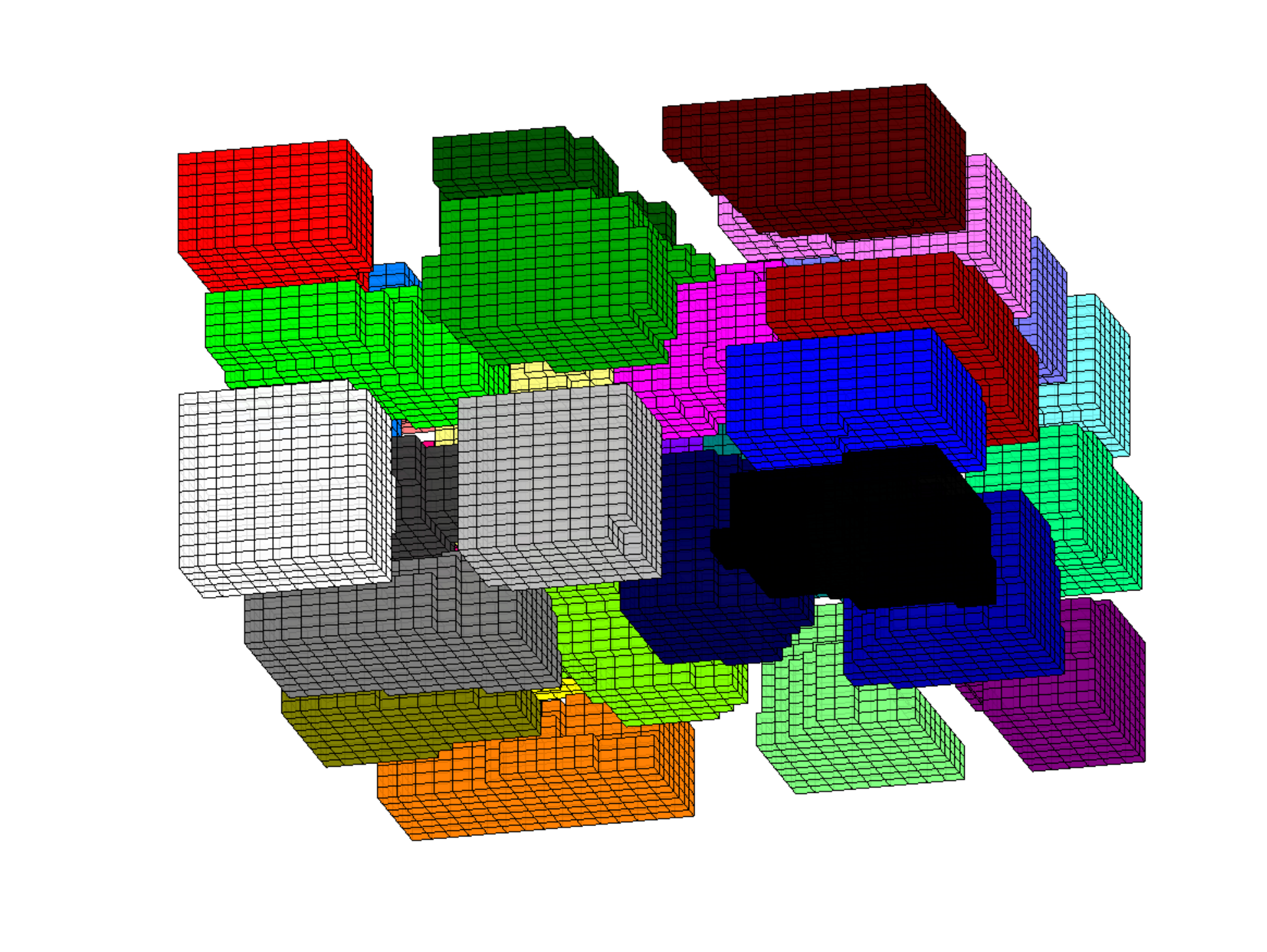}
\end{center}
\caption{Irregular partitioning of the domain from
Figure~\ref{fig:mixed_RT0_3D} into $32$ subdomains.}%
\label{fig:3D-metis}%
\end{figure}

We implemented the method in \textsc{Matlab}\ and studied
its convergence for problems with large variations in the permeability
coefficients $k$. In all experiments we used relative residual tolerance
$10^{-6}$ as the convergence criterion for the conjugate gradients. First, we
run a test with jumps in$~k$ aligned with substructure interfaces, see
Figure~\ref{fig:jumps}. For this problem we used stiffness scaling, which is
in case of the lowest-order Raviart-Thomas (RT0) elements equivalent\ to the
$\rho$-scaling. This also implies that the stiffness scaling works well for
irregular meshes (unlike for nodal elements). The conjugate gradients with the
BDDC\ preconditioner converged in~$15$ steps and the approximate condition
number computed from the L\'{a}nczos sequence in conjugate gradients was
$\kappa=4.046$; with $k=1$ the method converged in $14$ steps and
$\kappa=4.050$, see the rightmost column in Table~\ref{tab:layers}. In the
remaining experiments, we focused on problems with highly heterogeneous
coefficients, and we used the multiplicity scaling. Specifically, we simulated
flow in a porous media given by Model~2 of the 10th SPE Comparative Solution
Project~\cite{Christie-2001-SPE10}, which is publicly available on the
Internet\footnote{\url{http://www.sintef.no/Projectweb/GeoScale/Results/MsMFEM/SPE10/}
} and, in particular, we used a \textsc{Matlab} dataset described
in~\cite{Aarnes-2007-INF}. The dimensions of the full model are $1200\times
2200\times170$~(ft), and the distribution of the coefficients$~k$ is given
over a regular Cartesian grid with $60\times220\times85$ grid-blocks. We used
several layers and two 3D cutouts of the model for our numerical experiments.
For the experiments in 2D, we used layers 1, 20, 60 and 85 shown in
Figures~\ref{fig:layers_1_20}--\ref{fig:layers_60_85}. In the top layers~1 and
20 the permeability is relatively smooth, whereas the bottom layers~60, and 85
are fluvial and they are characterized by a spaghetti of narrow high-flow
channels. In all layers the permeabilities range over at least six orders of
magnitude. To drive a flow, we impose an injection (source) and a production
well (sink) in the lower-left and upper-right corners, respectively. The
discretization of each layer by the quadrilateral RT0 finite elements yields
$39,880$ degrees of freedom. The layers were partitioned into subdomains in
four ways: using two geometrically regular partitionings with the coarsening
ratios $H/h=30$ and $H/h=10$, and two irregular partitionings. The details of
the partitionings are summarized in Table~\ref{tab:partitioning}\ and
illustrated by Figures~\ref{fig:layers_1_20}--\ref{fig:layers_60_85}. For the
experiments in~3D, we used two domains consisting of $30\times30\times30$
elements extracted from layers~$1$--$30$ and $56$--$85$ of the SPE~10 problem
shown in Figure~\ref{fig:mixed_RT0_3D}. To drive a flow, we impose an
injection (source) and a production well (sink) in two distant corners of the
domain. The discretization by the hexahedral lowest-order Raviart-Thomas (RT0)
finite elements yields $110,700$ degrees of freedom. The domain was
partitioned into subdomains in two ways: using one geometrically regular
partitioning with the coarsening ratio $H/h=10$, and an irregular
partitioning. The details of the partitionings are summarized in
Table~\ref{tab:partitioning} and illustrated by Figure~\ref{fig:3D-metis}. All
irregular partitionings were obtained using \textsc{METIS~4.0}%
~\cite{Karypis-1998-MSP}, and in order to test the adaptive algorithm we did
not take into account the permeability coefficients.

It is interesting to note that the adaptive flux coarse basis functions
capture to some extent features of the solution on the finite element mesh,
and the quality of this approximation improves as the threshold$~\tau$ in
Algorithm~\ref{alg:adaptive} decreases. We illustrate this fact by relative
errors of solutions~$u_{0}$ and~$u^{\ast}$\ obtained in Steps~1 and~2 of
Algorithm~\ref{alg:two-level-nested} with respect to the exact
solution$~u_{\text{exact}}$ obtained by a direct solve of the full problem.
Specifically, the two relative errors are reported in tables as
\begin{equation}
\epsilon_{0}=\frac{\left\Vert u_{0}-u_{\text{exact}}\right\Vert }{\left\Vert
u_{\text{exact}}\right\Vert },\qquad\epsilon^{\ast}=\frac{\left\Vert u^{\ast
}-u_{\text{exact}}\right\Vert }{\left\Vert u_{\text{exact}}\right\Vert }.
\label{eq:epsilon}%
\end{equation}

We also compare the adaptive method with constraints inspired by
\textsl{Multiscale mixed finite element method (MsMFEM)}
cf.~\cite[Algorithm~2.5.2]{Efendiev-2009-MFE}
or~\cite[Section~3.2.1]{Aarnes-2006-HMM}.
In particular, instead of the local eigenvalue problems we solved local
Darcy's flow problems, that is local counterparts of
problem~(\ref{eq:problem-1})--(\ref{eq:problem-2}), with the source term
\[
f\left(  x\right)  =%
\begin{cases}
\;\;w_{i}, & \text{for }x\in\Omega^{i}\text{,}\\
-w^{j}, & \text{for }x\in\Omega^{j}\text{,}%
\end{cases}
\]
and zero flux boundary condition on $\partial\Omega^{i}\cap\partial\Omega^{j}%
$. The source distribution function is set to $w_{i}\left(  x\right)
=1/\left\vert \Omega^{i}\right\vert $ in all subdomains except those
containing a well, in which
\[
w_{i}\left(  x\right)  =\frac{f\left(  x\right)  }{\int_{\Omega^{i}}f\left(
\xi\right)  \,d\xi},
\]
to ensure a conservative approximation on the fine grid. In the numerical
experiments we then used the set of basic
constraints~(\ref{eq:starting_flux_constraint}) enriched by solving the above
problem and taking the values of flux degrees of freedom on $\partial
\Omega^{i}\cap\partial\Omega^{j}$ as additional constraints. Nevertheless, we
note that there are other more advanced solvers based on multiscale strategies
available in the literature see, e.g., Yang et al.~\cite{Yang-2018-TGP} or la
Cour Christensen et al.~\cite{laCourChristensen-2017-NMU}, and a thorough
comparison of the methods would be of independent interest.

\begin{figure}[ptbh]
\begin{center}
\includegraphics[width=5.5cm]{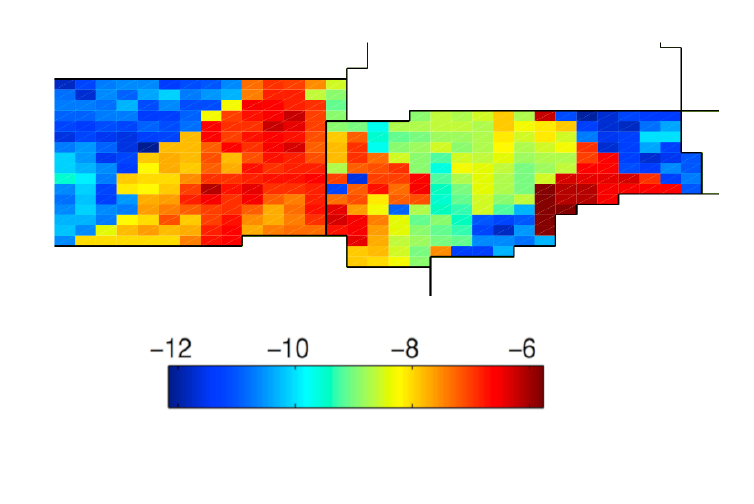}
\includegraphics[width=5.5cm]{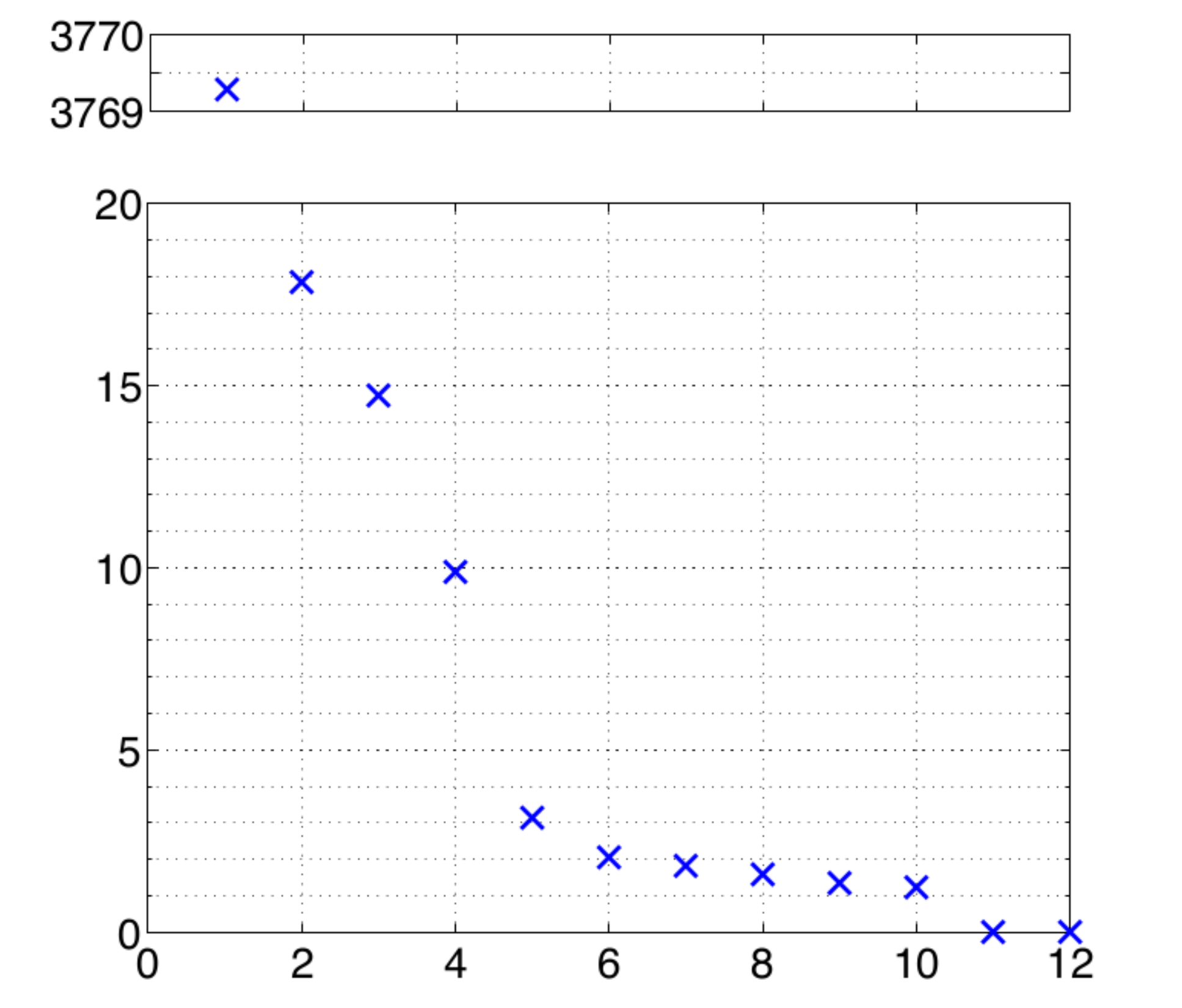}
\end{center}
\caption{The base~$10$ logarithm of the permeability~$k$ in the subdomains~$1$
and~$2$ from the layer~85 of the SPE~10 problem (Fig.~\ref{fig:layers_60_85},
right panel), and the first $12$ eigenvalues of the corresponding local
generalized eigenvalue problem~(\ref{eq:eigenvalue-problem-local}). Here
$\lambda_{i}=0$ for $i=11$ and $12$.}%
\label{fig:spe10_sub_1-2}%
\end{figure}

\begin{figure}[ptbh]
\begin{center}
\includegraphics[width=9cm]{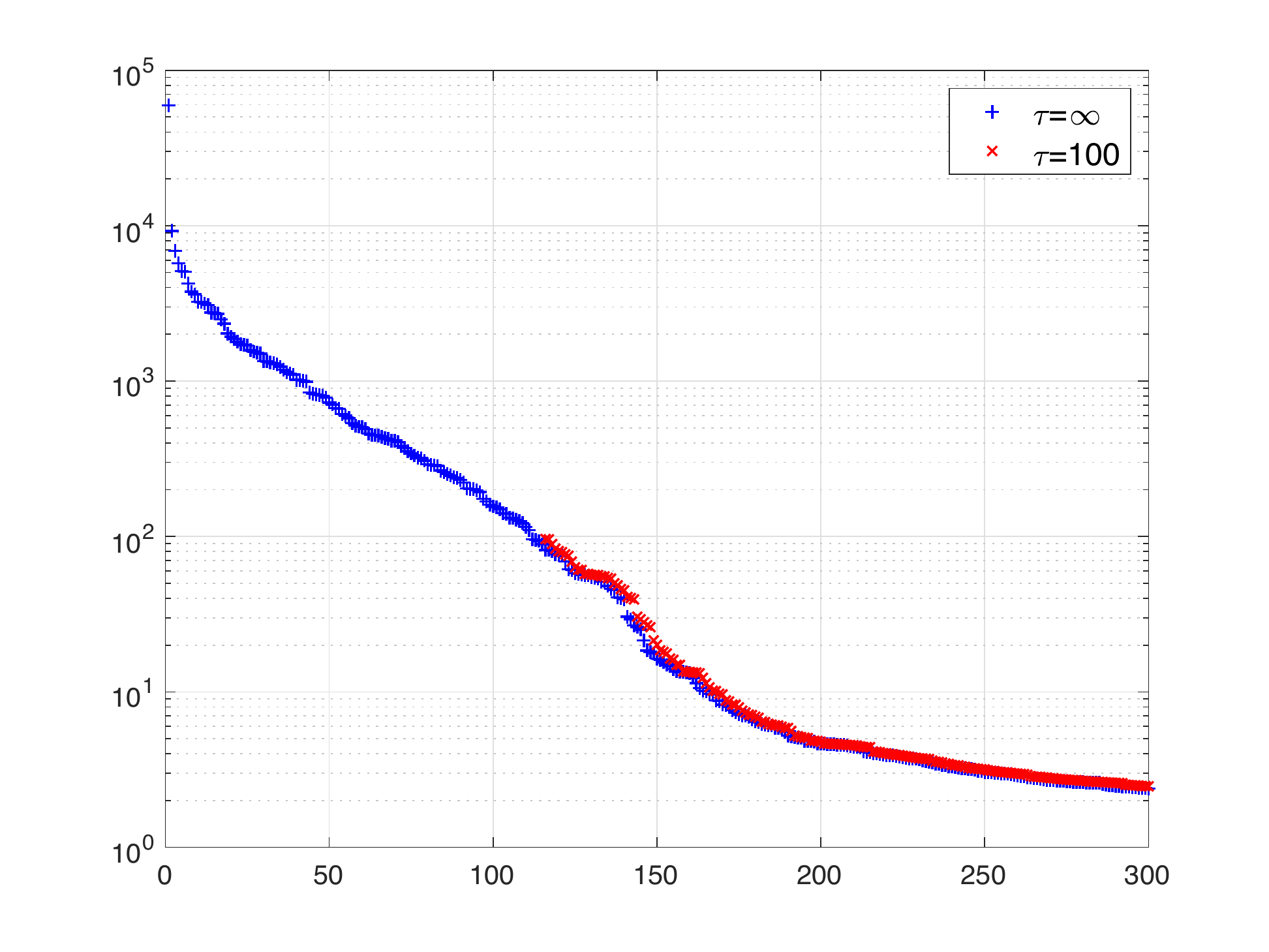}
\end{center}
\caption{The largest $300$ eigenvalues of the BDDC\ preconditioned operator
for the layer~$85$ of the SPE~10 problem (Fig.~\ref{fig:layers_60_85}, right
panel) without adaptivity ($\tau=\infty$) and for the adaptive BDDC with the
target condition number \thinspace$\tau=100$. }%
\label{fig:spe10_eig}%
\end{figure}

The results of numerical experiments in~2D are summarized in
Tables~\ref{tab:layers}--\ref{tab:layer85}. Table~\ref{tab:layers} shows
performance of the nonadaptive method for a homogeneous case with$~k=1$ and
the layers of the SPE~10 problem. It can be seen that for layers$~1$, and $20$
the convergence does not significantly depend on the partitioning and it is
also quite comparable to the homogeneous case with $k=1$. On the other hand,
for layers $60$ and $85$ the variations in coefficients aggravate convergence,
which is also quite sensitive to the partitioning. This holds, in particular,
for layer$~60$ which contains both regions that are highly heterogeneous and
relatively homogeneous. It can be also seen by comparing left and right
columns in Table~\ref{tab:layers} that increasing the number of subdomains
(that is decreasing the coarsening ratio$~H/h$) leads to higher condition
numbers and increase in iteration counts for both regular and irregular
partitionings. This is not the case in the standard theory of domain
decomposition methods, but here we suspect it can be attributed to the jumps
in coefficients and larger interfaces. The performance of the adaptive
algorithm is illustrated by Tables~\ref{tab:layer1}--\ref{tab:layer85}.
Table~\ref{tab:layer1} shows convergence for layer~$1$ with irregular
partitioning~A, and Table~\ref{tab:layer85} shows convergence for layer~$85$
with irregular partitioning~B. It can be seen that in all cases lower values
of the threshold$~\tau$ lead to fewer iterations, and the value of the
condition number indicator$~\widetilde{\omega}<\tau$ is in a good agreement
with$~\kappa$, which is the approximate condition number estimate obtained
from the L\'{a}nczos sequence in conjugate gradients. The adaptive constraints
also lead to more significant improvement in convergence than the multiscale
constraints. The problem for layer$~85$ is particularly interesting. From the
right panel in~Figure~\ref{fig:layers_60_85} we see that the coefficient jumps
have very large variations even on the interfaces, which can be seen in the
left panel of Figure~\ref{fig:spe10_sub_1-2}. The right panel displays the
eigenvalues of the corresponding eigenproblem: $\lambda_{1}\approx3769.5$ and
all other eigenvalues are less than$~20$. \ Figure~\ref{fig:spe10_eig} then
displays $300$ largest eigenvalues of the (global) BDDC\ preconditioned
operator without adaptivity and with adaptive BDDC\ and target condition
number $\tau=100$. We see that without adaptivity there is a single largest
eigenvalue: specifically $\lambda_{1}=59,492$ and $\lambda_{2}=9,258$. For the
adaptive BDDC\ with $\tau=100$ we get $\lambda_{1}=96.3$. Comparing this plot
with Table~\ref{tab:layer85} we see that the adaptive BDDC\ with $\tau=100$
introduces $115$ adaptive constraints, which corresponds to the number of the
largest eigenvalues removed from the spectrum of the BDDC preconditioned
operator. We also note that adding a single adaptive constraint reduces the
iteration count from $392$ to $347$, which corresponds to the large gap in the
spectrum of the operator without adaptivity. Setting $\tau$ to a lower value,
for example, $\tau=3$, roughly doubles the number of constraints and the
number of iterations is reduced to approximately~$10$. Also, the lower value
of$~\tau$ improve the approximation quality of the first two steps of
Algorithm~\ref{alg:two-level-nested} and, for example, with $\tau<3$ we get
the error $\epsilon^{\ast}<20\%$.

The results of numerical experiments in 3D are summarized in
Tables~\ref{tab:3D}--\ref{tab:spe10-3D-top}. It can be seen from
Table~\ref{tab:3D} that the numbers of iterations are significantly higher
than in~2D, and the convergence is slower for the fluvial bottom
layers~$56$--$85$ comparing with the relatively smooth top layers~$1$--$30$.
The increase in iterations becomes even more pronounced in the case of the
irregular partitioning also due to larger interfaces. The results of
experiments with the adaptive algorithm are summarized in
Tables~\ref{tab:spe10-3D}--\ref{tab:spe10-3D-top}. As in the~2D case, lower
values of the threshold$~\tau$ lead in all cases to fewer iterations, and the
values of$~\tau$, $\widetilde{\omega}$ and$~\kappa$ are in close agreement.
Again, the multiscale constraints provide only a slight improvement of
convergence. Table~\ref{tab:spe10-3D} shows convergence for layers$~1$--$30$.
It can be seen that despite higher condition number of the problem
corresponding to the irregular partitioning, the adaptive algorithm leads
allows to decrease the iteration counts for lower values of$~\tau$. As in 2D,
the first few adaptive constraints allow to decrease the number of iterations
by a fairly large amount: here adding $14$ constraints reduces the number of
iterations from the initial value~$1968$ to $1280$. However, for example
with$~\tau=10$ the number of iterations decreases to~$18$, however the number
of constraints grows rather significantly from~$335$ to~$1617$. Finally, the
values of$~\epsilon_{0}$ and$~\epsilon^{\ast}$ are quite larger compared to
the~2D experiments. Table~\ref{tab:spe10-3D-top} shows convergence for
layers$~56$--$85$ and the regular partitioning $H/h=10$, and the trends are
quite similar as in the previous case. That is, the adaptive algorithm allows
to control the convergence of conjugate gradients, but the number of adaptive
constraints is relatively high in particular for lower values of $\tau$. These
trends are in agreement with the qualitative observations made from
Figure~\ref{fig:spe10_eig}.

\begin{table}[ptbh]
\caption{Substructuring of the 2D and 3D problems: $N$ is the number of
subdomains, $n_{\Gamma}$ is the number of (flux) degrees of freedom on
interfaces, $n_{f}$ is the number of faces, and $n_{c}$ is the number of
(initial) coarse degrees of freedom.}%
\label{tab:partitioning}
\begin{center}%
\begin{tabular}
[c]{|c|r|r|r|r|}\hline
type of partitioning & $N$ & $n_{\Gamma}$ & $n_{f}$ & $n_{c}$\\\hline
\multicolumn{5}{|c|}{2D}\\\hline
regular ($H/h=30$) & 14 & 580 & 19 & 33\\
regular ($H/h=10$) & 132 & 2360 & 236 & 368\\
irregular A & 16 & 756 & 29 & 70\\
irregular B & 64 & 1746 & 152 & 315\\\hline
\multicolumn{5}{|c|}{3D}\\\hline
regular ($H/h=10$) & 27 & 5400 & 54 & 81\\
irregular & 32 & 7267 & 129 & 335\\\hline
\end{tabular}
\end{center}
\end{table}

\begin{table}[ptbh]
\caption{Convergence of the non-adaptive method for the homogeneous case
($k=1$) and the six layers of the SPE~10 problem. Here $\epsilon_{0}$ and
$\epsilon^{\ast}$ are the errors (in $\%$) defined by~(\ref{eq:epsilon}),
$\protect\widetilde{\omega}$ is the condition number indicator
from~(\ref{eq:indicator}), $it$~is the number of iterations for relative
residual tolerance $10^{-6}$, and $\kappa$~is the approximate condition number
computed from the L\'{a}nczos sequence in conjugate gradients.}%
\label{tab:layers}
\begin{center}%
\begin{tabular}
[c]{|c|r|r|r|r|r|r|r|r|}\hline
\multirow{2}{*}{layer} & \multicolumn{2}{|c|}{$H/h=30$} &
\multicolumn{2}{|c|}{$H/h=10$} & \multicolumn{2}{|c|}{irregular~A} &
\multicolumn{2}{|c|}{irregular~B}\\\cline{2-9}
& $it$ & $\kappa$ & $it$ & $\kappa$ & $it$ & $\kappa$ & $it$ & $\kappa
$\\\hline
($k=1$) & 11 & 2.790 & 14 & 3.980 & 12 & 3.151 & 14 & 4.050\\
1 & 15 & 8.879 & 22 & 9.491 & 17 & 6.714 & 19 & 11.197\\
20 & 14 & 5.749 & 19 & 6.926 & 15 & 6.524 & 18 & 6.429\\
60 & 162 & 4564.1 & 513 & 26,359.3 & 244 & 11,272.6 & 292 & 7301.7\\
85 & 183 & 9310.7 & 446 & 24,492.8 & 208 & 7170.4 & 392 & 58,931.7\\\hline
\end{tabular}
\end{center}
\end{table}

\begin{table}[ptbh]
\caption{Convergence of the adaptive method for layer~$1$ of the SPE~10
problem with the irregular partitioning A. Here $\tau$ is the target condition
number from Algorithm~\ref{alg:adaptive}, $\epsilon_{0}$~and~$\epsilon^{\ast}$
are the errors (in~$\%$) defined by~(\ref{eq:epsilon}),
$\protect\widetilde{\omega}$ is the condition number indicator
from~(\ref{eq:indicator}), $it$~is the number of iterations for relative
residual tolerance $10^{-6}$, and $\kappa$~is the approximate condition number
computed from the L\'{a}nczos sequence in conjugate gradients. With
$\tau=\infty$ no adaptive constraints were used, and (ms) indicates use of the
multiscale constraints.}%
\label{tab:layer1}
\begin{center}%
\begin{tabular}
[c]{|r|r|r|r|r|r|r|}\hline
$\tau$ & $\epsilon_{0} \,\, [\%]$ & $\epsilon^{*} \,\, [\%]$ &
$\widetilde{\omega}$ & $n_{c}$ & $it$ & $\kappa$\\\hline
$\infty$ & 73.21 & 30.55 & 13.586 & 70 & 17 & 6.714\\\hline
(ms) & 72.55 & 27.15 & -na-$\,\,$ & 121 & 15 & 5.998\\\hline
10 & 71.86 & 29.11 & 8.404 & 73 & 16 & 6.231\\
5 & 70.77 & 18.89 & 4.765 & 81 & 13 & 5.517\\
3 & 69.19 & 11.64 & 2.997 & 104 & 11 & 2.842\\
2 & 69.23 & 9.42 & 1.970 & 153 & 8 & 1.915\\\hline
\end{tabular}
\end{center}
\end{table}

\begin{table}[ptbh]
\caption{Convergence of the adaptive method for layer~$85$ with the irregular
partitioning~B. }%
\label{tab:layer85}
\begin{center}%
\begin{tabular}
[c]{|r|r|r|r|r|r|r|}\hline
$\tau$ & $\epsilon_{0} \,\, [\%]$ & $\epsilon^{*} \,\, [\%]$ &
$\widetilde{\omega}$ & $n_{c}$ & $it$ & $\kappa$\\\hline
$\infty$ & 69.34 & 42.11 & 59,491.702 & 315 & 392 & 58,931.700\\\hline
(ms) & 68.76 & 39.00 & -na-$\quad$ & 494 & 297 & 8931.930\\\hline
10,000 & 69.34 & 42.11 & 9275.614 & 316 & 347 & 9170.830\\
1000 & 68.03 & 40.29 & 898.754 & 360 & 152 & 836.227\\
100 & 67.21 & 38.38 & 98.117 & 430 & 54 & 95.439\\
10 & 66.16 & 35.68 & 9.885 & 489 & 19 & 9.672\\
5 & 63.07 & 31.85 & 4.888 & 536 & 13 & 4.836\\
3 & 56.19 & 18.87 & 2.988 & 614 & 10 & 3.010\\
2 & 53.44 & 14.87 & 1.997 & 743 & 7 & 1.879\\\hline
\end{tabular}
\end{center}
\end{table}

\begin{table}[ptbh]
\caption{Convergence of the non-adaptive method for the homogeneous case
($k=1$) and the two 3D cutouts of the SPE~10 problem from
Figure~\ref{fig:mixed_RT0_3D}. The headings are same as in
Table~\ref{tab:layers}.}%
\label{tab:3D}
\begin{center}%
\begin{tabular}
[c]{|c|r|r|r|r|}\hline
\multirow{2}{*}{layer} & \multicolumn{2}{|c|}{$H/h=10$} &
\multicolumn{2}{|c|}{irregular part.}\\\cline{2-5}
& $it$ & $\kappa$ & $it$ & $\kappa$\\\hline
$(k=1)$ & 25 & 17.099 & 35 & 22.029\\
1--30 & 779 & 49,075.600 & 1968 & $1.096\times10^{6}$\\
56--85 & 3762 & $2.576\times10^{6}$ & 5277 & $3.676\times10^{6}$\\\hline
\end{tabular}
\end{center}
\end{table}

\begin{table}[ptbh]
\caption{Convergence of the adaptive method for layers~1--30 of the SPE~10
with irregular partitioning. }%
\label{tab:spe10-3D}
\begin{center}%
\begin{tabular}
[c]{|r|r|r|r|r|r|r|}\hline
$\tau$ & $\epsilon_{0} \,\, [\%]$ & $\epsilon^{*} \,\, [\%]$ &
$\widetilde{\omega}$ & $n_{c}$ & $it$ & $\kappa$\\\hline
$\infty$ & 98.41 & 86.05 & $1.191\times10^{6}$ & 335 & 1968 & $1.096\times
10^{6}$\\\hline
(ms) & 98.29 & 86.05 & -na-$\quad\,$ & 571 & 1943 & $1.079\times10^{6}%
$\\\hline
100,000 & 98.39 & 85.95 & 94,328.862 & 349 & 1280 & 92,307.000\\
10,000 & 98.14 & 84.27 & 9862.559 & 514 & 514 & 10,512.200\\
1000 & 97.41 & 82.73 & 995.230 & 989 & 175 & 1014.150\\
100 & 92.93 & 72.63 & 97.989 & 1331 & 60 & 108.673\\
10 & 87.47 & 66.66 & 9.993 & 1617 & 18 & 11.711\\
5 & 85.82 & 65.46 & 4.985 & 1898 & 13 & 6.069\\
3 & 82.90 & 63.22 & 2.997 & 2331 & 9 & 3.007\\
2 & 81.62 & 62.81 & 2.000 & 2997 & 6 & 1.930\\\hline
\end{tabular}
\end{center}
\end{table}

\begin{table}[ptbh]
\caption{Convergence of the adaptive method for layers~56--85 of the SPE~10
with the regular partitioning. }%
\label{tab:spe10-3D-top}
\begin{center}%
\begin{tabular}
[c]{|r|r|r|r|r|r|r|}\hline
$\tau$ & $\epsilon_{0} \,\, [\%]$ & $\epsilon^{*} \,\, [\%]$ &
$\widetilde{\omega}$ & $n_{c}$ & $it$ & $\kappa$\\\hline
$\infty$ & 99.31 & 66.24 & $3.435\times10^{6}$ & 81 & 3762 & $2.576\times
10^{6}$\\\hline
(ms) & 99.30 & 66.24 & -na-$\quad\,$ & 99 & 3735 & $2.566\times10^{6}$\\\hline
100,000 & 99.33 & 65.78 & 95,129.959 & 122 & 1267 & 93,040.500\\
10,000 & 98.63 & 62.19 & 9834.429 & 188 & 498 & 9487.840\\
1000 & 98.25 & 63.81 & 990.920 & 373 & 183 & 1200.070\\
100 & 97.26 & 64.08 & 99.793 & 766 & 59 & 124.896\\
10 & 91.79 & 49.39 & 9.965 & 1154 & 17 & 9.506\\
5 & 88.88 & 43.77 & 4.991 & 1342 & 13 & 5.545\\
3 & 86.52 & 41.75 & 2.997 & 1610 & 9 & 3.205\\
2 & 85.18 & 41.05 & 2.000 & 1960 & 6 & 1.918\\\hline
\end{tabular}
\end{center}
\end{table}

\section{Conclusion}

\label{sec:conclusion}We studied a method for solution of single-phase flow in
heterogeneous porous media. We have, in particular, shown that the idea of
adaptive BDDC, previously used for elliptic problems, can be also applied in
the context of the BDDC method for mixed finite element discretizations using
the lowest-order Raviart-Thomas finite elements, and that the adaptive method
works well with the usual types of scaling used in substructuring. We
illustrated that the resulting algorithm can be successfully applied for
adaptive selection of the coarse flux degrees of freedom using several
examples corresponding to the SPE~10 benchmark model. The effect of the
adaptive construction of the flux coarse basis functions is twofold. First,
the first two steps of the BDDC method provide some approximation properties
with respect to the exact solution of the full problem, in particular in 2D.
Second, the coarse problem provides a better preconditioner for conjugate
gradients used in the third step. We also compared the adaptive constraints
with constraints inspired by multiscale mixed finite element method, and we
found that the adaptive constraints outperform the multiscale constraints.

Next, we experimented with different partitionings of the domains into
substructures. While the adaptive method is able to overcome these issues in
many cases, it is evident that a suitable partitioning makes the adaptive
method more efficient. We note that development of optimal partitioning
strategies is an open problem cf.,
e.g.,~\cite{Aarnes-2008-MMM,Vecharynski-2014-GPU}. However, our experiments
indicate that if it is not possible to find a suitable partitioning, the best
strategy is to simply minimize the size of interfaces, which may be achieved
by a simple geometric partitioning, see also~\cite{Hanek-2017-EII}.

\newpage


\bibliographystyle{siam}
\bibliography{bddc_spe10}

\end{document}